\documentclass[sigconf]{acmart}
\AtBeginDocument{%
  }

\setcopyright{acmlicensed}
\copyrightyear{2026}
\acmYear{2026}
\acmDOI{10.1145/3744255.3811725}
\setcctype{by}
\acmConference[E-Energy '26]{The 17th ACM International Conference on Future and Sustainable Energy Systems}{June 22--25, 2026}{Banff, AB, Canada}
\acmBooktitle{The 17th ACM International Conference on Future and Sustainable Energy Systems (E-Energy '26), June 22--25, 2026, Banff, AB, Canada}
\acmISBN{979-8-4007-2011-6/2026/06}

\usepackage{amsmath}
\usepackage{graphicx}
\usepackage{subcaption}

\allowdisplaybreaks

\begin{document}

\title{Robust Capacity Expansion under Wildfire Ignition Risk and High Renewable Penetration}

\author{Tomás Tapia}
\affiliation{%
  \institution{Johns Hopkins University}
  \city{Baltimore}
  \state{Maryland}
  \country{USA}
}

\author{Ryan Piansky}
\affiliation{%
  \institution{Georgia Institute of Technology}
  \city{Atlanta}
  \state{Georgia}
  \country{USA}
}

\author{Yury Dvorkin}
\affiliation{%
  \institution{Johns Hopkins University}
  \city{Baltimore}
  \state{Maryland}
  \country{USA}
}

\author{Jean-Paul Watson}
\affiliation{%
  \institution{Lawrence Livermore National Laboratory}
  \city{Livermore}
  \state{California}
  \country{USA}
}

\renewcommand{\shortauthors}{T. Tapia, R. Piansky, Y. Dvorkin, and J-P. Watson}

\begin{abstract}
    In power systems, the risk of wildfire ignition has increased significantly in recent years. The impact and severity of these events on energy dispatch, as well as their societal ramifications, make wildfire prevention critical for power system planning and operation. A common intervention by system operators is to de-energize transmission lines to mitigate the risk of fire caused by equipment failures. With the growing integration of variable renewable generation, managing and preparing the system to de-energization under wildfire risk has become even more challenging. In this context, mitigation decisions such as installing battery energy storage systems and undergrounding transmission lines can reduce the risk and adverse effects associated with de-energization and renewable generation variability. This paper presents a robust optimization model to determine the optimal location of battery storage and undergrounding of transmission line investment, utilizing representative weeks and uncertainty sets to capture the temporal relationship of uncertain variables. Specifically, this paper addresses: (i) the worst-case realization of ignition risk leading to the de-energization of transmission lines, combined with the worst-case realization of renewable energy availability, and (ii) the optimal investment decisions for energy storage capacity and undergrounding of transmission lines that are exposed to ignition risk. The proposed model is formulated as a mixed-integer linear programming (MILP) problem, employing duality theory and binary decomposition to address nonlinearities, and is solved using a column-and-constraint generation algorithm. The proposed framework is evaluated on a model of the San Diego power system, demonstrating its practical effectiveness in improving the resilience to wildfire risk.
\end{abstract}




\maketitle

\section{Introduction}
\subsection{Motivation and Scope}
Due to the continued evolution of power systems, capacity expansion planning (CEP) is a key process for governmental energy commissions and utilities to ensure resource adequacy, reliability, and cost-efficient energy supply. In many jurisdictions, CEP is also mandated by law as part of regulatory compliance. CEP finds the least-cost mix of new investments in generation, storage, and transmission infrastructure of a power system over a given planning horizon, based on an estimate of future system requirements such as demand growth, renewable integration targets, reliability standards, and policy constraints. With the increasing integration of variable sources and the growing occurrence of of weather-driven events, this process has become challenging due to the significant levels and sources of uncertainty introduced to the system \cite{musselmanclimate}. Several states, including California, Texas, and New York, are increasingly adopting modeling approaches to explicitly account for these uncertainties in their long-term planning processes. 

Underestimating power system variability factors, such as fluctuating renewable generation and load demand, generator outages, fuel supply disruptions, transmission constraints, and extreme weather events, can lead to significant contingencies that disrupt normal grid operations, potentially resulting in component disconnections and compromising the security of power supply. Given this, developing efficient tools that capture and cope with the inherent uncertainty has become a critical task for energy planners. In this context, wildfires represent a major threat to the secure operation of power systems. During recent years, California, Canada, Australia, Portugal, and Chile have suffered from wildfires that have lead to supply interruption \cite{tapia2021robust}. The treatment of wildfire events in power systems has been approached in different directions: (i) wildfire mitigation, focused on the pre-disaster stage that considers all forecasting and infrastructure investment; (ii) wildfire response, focused on emergency response actions; and (iii) disaster restoration, focused on infrastructure repair and energy rehabilitation strategies. Furthermore, as the adoption of variable renewable energy sources grows, strengthening their resilience to wildfires becomes more challenging due to the increased uncertainty these systems face. Unlike dispatchable fossil plants, wind and solar plants rely on weather and long lines, so wildfire smoke, heat, and line outages can simultaneously cut production and access, reducing resilience.

In the first direction, wildfires caused by the infrastructure of the power system (e.g. sparks, faults, etc.) are not atypical due to aging components and severe climate conditions that increase the potential failure, ignition, and fire spread \cite{piansky2025optimizing}. To mitigate this threat, utilities can temporarily turn off power or de-energize specific areas under high wildfire ignition risk to prevent infrastructure-related wildfire ignitions, the so-called \emph{Public Safety Power Shutoff} (PSPS) events \cite{cpuc_psps}. Consequently, system operators must carefully balance operating costs, manage component de-energization, and reduce wildfire ignition risk. In this context, accurately quantifying wildfire ignition risk can support more informed and effective component de-energization planning, particularly for transmission lines.

With ignition risk map information, system operators can better anticipate periods and locations of elevated wildfire risk and estimate periods where transmission lines will be de-energized to avoid operations during times of high ignition threat. Integrating this temporal risk information into CEP offers an opportunity to co-optimize both wildfire mitigation measures, such as preemptive transmission de-energization, and power system investments. 

This paper uses robust optimization to consider the uncertainty sets for wildfire ignition risk and renewable energy variability in CEP, ensuring that the investment takes into account the variability between ignition risk and variable generation availability forecasts under an uncertainty budget, potentially reducing the probability of catastrophic events while maintaining reliable power delivery. The resulting model is a mixed-integer linear program (MILP) that captures both the discrete nature of investment, operational decisions, and inter-temporal constraints.

\subsection{Literature review}
CEP problems have been widely studied due to their significant role in guiding investment and ensuring reliability. In recent years, various approaches have been explored to address uncertainty factors inherent to these problems, including stochastic optimization \cite{piansky2025optimizing}, robust optimization \cite{verastegui2019adaptive, tapia2021robust}, and distributionally robust optimization \cite{hong2022data}. Each of these methodologies employs different techniques to quantify uncertainty arising from variable factors such as weather conditions, renewable resource availability, and energy demand. For example, scenario-based methods \cite{musselmanclimate, piansky2025optimizing}, and uncertainty sets \cite{verastegui2019adaptive, tapia2021robust}, are commonly used to represent possible realizations of these uncertainties. 

The increasing integration of renewable energy sources into the grid, coupled with the increasing frequency of weather-driven events, has heightened the importance of accurately modeling these uncertainties to assess and enhance system resilience. As power systems continue to evolve towards a larger integration of variable generation energy, ensuring their robustness against weather-driven events becomes an even greater challenge.

A key challenge arises when line de-energization due to wildfire risk overlaps with low renewable availability, leading to suboptimal use of storage or investment decisions. CEP must also address supply and demand variability, growing flexibility needs, and the proper timing of mitigation actions. Failing to account for these factors can result in underinvestment or resource shortfalls. Representing wildfire risk and renewable variability helps identify when and where the system is most stressed, enabling more resilient planning. Our model addresses these challenges by incorporating renewable availability and wildfire risk into two different uncertainty sets.

In particular, events such as wildfires have received significant attention in recent research due to their complex evolution and potential impacts on multiple systems, making them highly relevant for consideration. Similarly to other resilience analyses, studies focusing on wildfires can be categorized into three main research directions: operational response, system planning, and restoration. The operational response involves actions such as transmission line adjustment, generation redispatch, and resource allocation and management (e.g. \cite{mohagheghi2015optimal, abdelmalak2022enhancing,piansky2024quantifying,Taylor2023Managing}]) to mitigate immediate operational impacts. In contrast, system planning focuses on long-term risk mitigation strategies (e.g. \cite{vazquez2022wildfire,taylor2023robust}). Finally, restoration focuses on recovering system functionality after wildfire events (e.g. \cite{rhodes2023co}). Within the context of wildfire mitigation, capacity expansion planning often incorporates various risk metrics with the objective of minimizing unserved energy. Recent studies have addressed wildfire ignition risk by defining and comparing different metrics, utilizing both threshold-based and optimization-based methods \cite{piansky2024quantifying}. Authors in \cite{kody2022optimizing, bayani2023resilient} analyze the co-optimization of line de-energization and transmission investments, while \cite{piansky2025optimizing} proposes a scenario-based optimization framework for battery investments and undergrounding transmission lines. In addition, robust optimization frameworks have been developed to analyze transmission switching strategies under worst-case wildfire spread scenarios and varying renewable energy availability \cite{tapia2021robust}. 

Authors in \cite{tapia2021robust} constructed uncertainty sets related to wildfire risk, but their focus was on components failure driven by fire-spread propagation across the power network. In contrast to this approach, our uncertainty set representation focuses on wildfire ignition risk to support mitigation decisions.

This paper addresses a gap in the literature by focusing on power system planning under worst-case wildfire ignition risk and renewable availability scenarios. We consider the de-energization of transmission lines as a mitigation measure against wildfire impacts and construct tailored uncertainty sets that capture both ignition risk and renewable variability. This representation is well-suited for planning problems because it captures the primary system-level effects of wildfire risk through a tractable model, while enabling more informed and resilient capacity expansion decisions under extreme but plausible conditions. Nevertheless, it may be conservative, as the selected line-risk metric and de-energization threshold reduce the continuous nature of actual ignition risk to a binary approximation.

\subsection{Contributions}
This paper develops a two stage algorithm based on robust optimization to solve CEP considering impacts of the worst-case transmission line de-energization and renewable availability realization. The main contributions are:
\begin{itemize}
    \item A robust framework for wildfire ignition-aware capacity expansion planning is presented for electricity grids with high renewable penetration. The proposed method incorporates ignition risk metrics associated with each transmission line. The model jointly balances transmission and storage investment decisions while preemptively accounting the de-energization of transmission lines and worst-case realization renewable generation availability.
    \item The proposed model uses two traditional uncertainty sets for transmission lines de-energization, and renewable trajectories. Thus, enables utilities and power system planners to calibrate risk tolerance corresponding budget of uncertainty. These budgets define the allowable deviation from nominal conditions within each uncertainty set, controlling the level of conservatism in the planning solution. Thus, out approach supports the identification of an optimal investment plan alongside the adaptive response on generator and storage dispatch actions.
    \item Extensive computational experiments on a test case based on the San Diego system with up to 24 buses the practical applicability of the wildfire ignition-aware framework proposed in this paper. With this, several cases of budget values are presented to analyze the impacts of uncertainty realization, transmission line de-energization, and power system expansion decisions.
\end{itemize}

\section{Wildfire ignition risk maps and ignition risk representation} \label{section:ignition_maps}
There are several data sets focused on quantifying wildfire potential, typically using a combination of surface meteorological measurements and environmental factors such as humidity, fuel type, wind speed, and topography \cite{usgs_wfpi_2025,pyregence}. These models generally pursue two main objectives: (i) quantifying the risk of wildfire ignition \cite{piansky2024quantifying, bayani2022quantifying}, and (ii) simulating the evolution of wildfire spread \cite{allaire2021emulation, tapia2021robust}. 

For our purposes, quantifying the risk of wildfire ignition is critical to implementing mitigation measures in power system operations. As noted in \cite{piansky2024quantifying}, Wildland Fire Potential Index (WFPI) Maps provided by the USGS \cite{usgs_wfpi_2025} are a strong candidate for assessing wildfire ignition risk associated with component failure of the power system. This is due to the WFPI’s daily forecast updates and its detailed spatial granularity. The maps present leveled risk information, i.e. categorized in a way that enables rapid interpretation and supports timely decision-making, for the entire United States, based on an algorithm that incorporates estimates of live-to-dead fuel ratios, wind speed, dry bulb temperature, and rainfall among other factors \cite{usgs_wfpi_2025}.

\begin{figure}[h]
    \centering
    \includegraphics[width=0.96\linewidth, trim={0, 0, 0, 60}, clip]{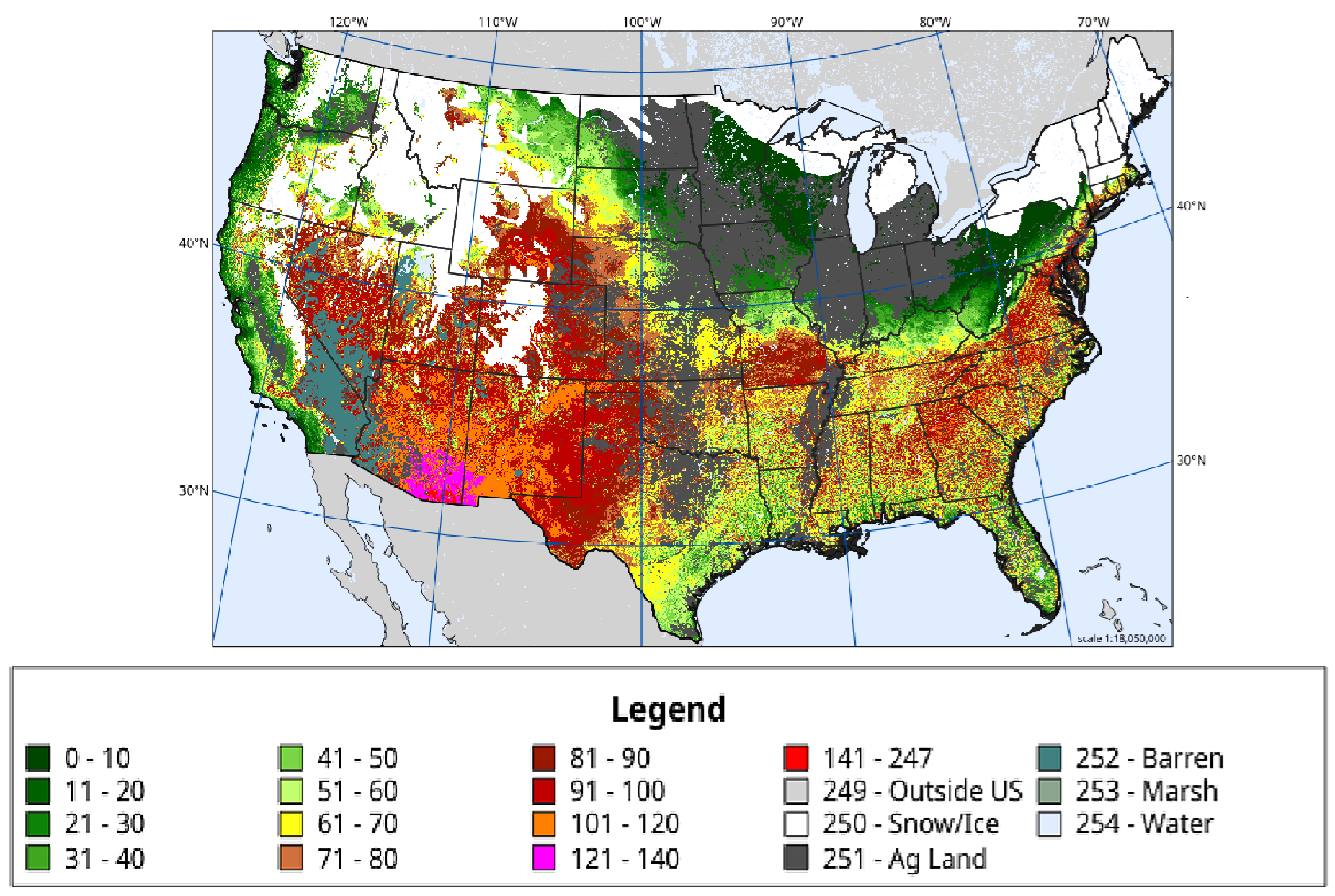}
    \caption{Wildland Fire Potential Index (WFPI) for the United States, March 1, 2025 \cite{usgs_wfpi_2025}.}
    \label{Fig:Map}
\end{figure}

To effectively incorporate wildfire ignition risk into CEP, we consider the use of representative weeks with uncertainty sets to represent the uncertainty and variability of the risk measures, instead of relying on scenario-based which can become computationally burdensome. The inclusion of key parameters from a representative week allows us to capture the temporal correlation while maintaining the problem to be tractable for large-scale network instances. Figure \ref{Fig:Map} presents an example of a WFPI wildfire-ignition risk map, which contains the risk measures that we will use to construct the risk uncertainty set.

Representative weeks of ignition risk are constructed by applying statistical analysis to 12 months of the WFPI Maps data, capturing both the estimated mean and deviation of risk values throughout the period. Specifically, we use the $k$-medoids algorithm to identify the set of representative weeks that best capture the diversity of risk profiles observed in the original data. Each representative week is assigned a corresponding weight, reflecting its frequency or importance in historical information. The $k$-medoids method is a clustering algorithm related to k-means, but instead of using the mean or average of cluster points as the center, it uses actual data points (so-called \textit{medoids}) as cluster centers. Its goal is to partition the data set into $k$ clusters, minimizing the sum of distances between the points and their respective medoids. Alternatively, the representative set can be obtained using a combination of $k$-means and $k$-medoids, as proposed in \cite{li2022representative}.

\begin{figure}[htbp]
    \centering
    \begin{subcaptionbox}{Week cluster 1, weight $=0.5$\label{fig:subfig1}}[.3\linewidth]
        {\includegraphics[width=\linewidth, trim={0, 0, 125, 0}, clip]{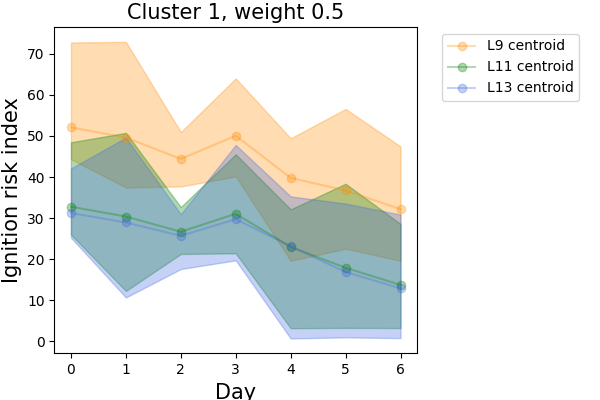}}
    \end{subcaptionbox}
    \hfill
    \begin{subcaptionbox}{Week cluster 2, , weight $=0.417$\label{fig:subfig2}}[.3\linewidth]
        {\includegraphics[width=\linewidth, trim={0, 0, 125, 0}, clip]{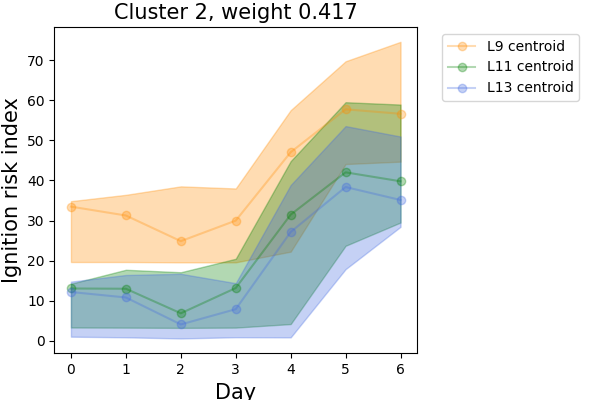}}
    \end{subcaptionbox}
    \hfill
    \begin{subcaptionbox}{Week cluster 3, weight $=0.083$\label{fig:subfig3}}[.3\linewidth]{\includegraphics[width=\linewidth, trim={0, 0, 125, 0}, clip]{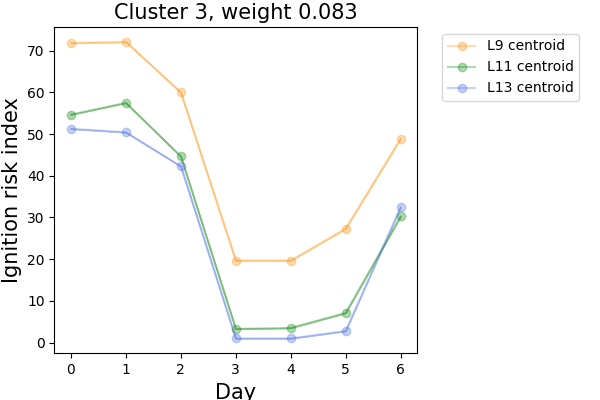}}
    \end{subcaptionbox}
    \caption{Illustration of three representative weeks. Each color represent a group of historical weeks, the line is the medoids trajectory and the colored area is the feasible risk variation.}
    \label{fig:representative}
\end{figure}

For example, Figure \ref{fig:representative} shows the representative weeks for three transmission lines taking information from three months (from March to May 2025) of WFPI maps. Thus, this approach provides a tractable yet informative foundation for subsequent capacity expansion and operational analysis, balancing computational efficiency with fidelity to observed risk patterns.

\section{Mathematical formulation for the ignition-Aware Robust CEP model}
In power systems, CEP models determine the optimal investment of new generation, storage, and transmission infrastructure over a planning horizon to satisfy future energy demands while maintaining reliable system operation. In the context of wildfire ignition risk, the objective changes to ensure not only reliable and cost-effective energy supply but minimizing the risk of starting a fire. The key challenge lies in maintaining the balance between implementing risk mitigation measures without significantly increasing capital and operating costs or compromising system reliability.

CEP models are typically formulated as mixed-integer optimization problems, accounting for investment costs, operational efficiency, exposure to various sources of uncertainty, and, in this case, exposure to wildfire-related risks. 

In this context, the proposed model adopts a robust optimization framework to protect the system against the worst-case realization of wildfire ignition risk and variable generation availability. It determines the investment costs for battery storage and line undergrounding, while ensuring both economic efficiency and reliable system operation. We summarize our notation in Appendix A. The following presents the proposed robust formulation:
\begin{subequations} \label{p0:invest_problem}
\begin{align} 
    \min_{X,z^{\rm UG}} ~& \sum_{l} C_{l}^{\rm UG} z_{l}^{UG} + \sum_{s} C_{s}^{X}~X_{s}^{\rm max} \nonumber \\
    & + \max_{z \in \mathcal{U}_{r}(z^{\rm UG}),~ \alpha \in  \mathcal{U}_{\alpha}}~\min_{p,f,e,d,c} f(\cdot), \label{p0:invest_objective}\\
    \text{s.t.}\quad &0 \leq X_{s}^{\rm max} \leq X_{s}^{\rm sup}, \label{p0:invest_limits} \\
    & z_{l}^{\rm UG} \in \{0,1\},    
\end{align}
\end{subequations}
$X_{s}^{\rm max}$ represents the storage capacity investment decision for candidate project $s \in S^{C}$ and $z_{l}^{UG}$ is the underground line investment decision for line $l$. Eq. \eqref{p0:invest_objective} corresponds to the capacity expansion objective function. It aims to minimize the inversion on underground lines and battery storage installation while minimizing the worst-case operational cost represented by $\max\min f(\cdot)$. The lower-level, $\min f(\cdot)$, corresponds to the minimization of the operating cost under a particular realization of the uncertainty. Eq. \eqref{p0:invest_limits} corresponds to the maximum capacity limit for candidate energy storage projects. The lower-level optimization problem is described below:
\begin{subequations} \label{p1:economic_dispatch}
\begin{align} 
    \min_{p,f,e,d,c}~& \sum_{w} w_{w} \sum_{t,d}\Big[ \sum_{g} C_{g} p_{g}^{t,d,w}  \nonumber\\
    & + \sum_{s} C_{s}^{\rm dis} d_s^{t,d,w} + \sum_{b} C_{b}^{LS} p_{b}^{t,d,w} \Big] \label{p1:objective_funtion}\\
    \text{s.t.}\quad& p_{g}^{t,d,w},~f_{l}^{t,d,w},~p_{b}^{t,d,w},~e_{s}^{t,d,w} \geq 0, \label{p1:variables}\\
    (\underline{\Delta}_{g}^{t,d,w},\overline{\Delta}_{g}^{t,d,w}):\quad& p_{g}^{\rm min} \leq p_{g}^{t,d,w} \leq p_{g}^{\rm max}, \label{p1:gen_limts}\\
    (\underline{\Delta}_{r}^{t,d,w},\overline{\Delta}_{r}^{t,d,w}):\quad & 0 \leq p_{r}^{t,d,w} \leq \alpha_{r}^{t,d,w} p_{r}^{\rm max}, \label{p1:gen_renewable_limts}\\
    (\pi_{b}^{t,d,w}):\quad& \sum_{g \in G_b} p_{g}^{t,d,w} + \sum_{r \in R_b} p_{r}^{t,d,w} \nonumber\\
    & + \sum_{s \in S_b} (d_{s}^{t,d,w} - c_{s}^{t,d,w}) \nonumber\\
    & + \sum_{l|b^{\rm in}(l) = b} f_{l}^{t,d,w} - \sum_{l|b^{\rm out}(l) = b} f_{l}^{t,d,w} \nonumber \\
    & = d_{b}^{t,d,w} - p_{b}^{t,d,w},\label{p1:energy_balance}\\
    (\underline{\Delta}_{l}^{t,d,w},\overline{\Delta}_{l}^{t,d,w}):\quad&  |f_{l}^{t,d,w}| \leq \overline{f}_{l} (1-z_{l}^{d,w}), \label{p1:line_limits}\\
    (\underline{\Theta}_{l}^{t,d,w}, \overline{\Theta}_{l}^{t,d,w}):\quad & |f_{l}^{t,d,w} - B_l(\theta_{b^{\rm in}(l)}^{t,d,w} - \theta_{b^{\rm out}(l)}^{t,d,w})| \leq M z_{l}^{d,w}, \label{p1:line_angles}\\
    (\underline{\Delta}_{s}^{t,d,w},\overline{\Delta}_{s}^{t,d,w}):\quad & 0 \leq e_{s}^{t,d,w} \leq X_{s}^{\rm max}, \label{p1:battery_limits}\\
    (\mu_{s}^{t,d,w}):\quad & e_{s}^{t+1,d,w} = \nonumber \\
    & e_{s}^{t,d,w} + \eta c_{s}^{t+1,d,w} - \frac{1}{\eta} d_{s}^{t+1,d,w},\label{p1:battery_dynamics}\\
    (\underline{\Theta}_{s}^{t,d,w},\overline{\Theta}_{s}^{t,d,w}):\quad  & 0 \leq d_{s}^{t,d,w} \leq P_{s}, \label{p1:battery_dis_limit}\\
    (\underline{\Phi}_{s}^{t,d,w},\overline{\Phi}_{s}^{t,d,w}):\quad  & 0 \leq c_{s}^{t,d,w} \leq P_{s}, \label{p1:battery_ch_limit}\\
    & e_{s}^{1,1,1} =  E_{s}^{\rm init}, \quad e_{s}^{T,D,W} =  E_{s}^{\rm end}. \label{p1:battery_initial}
\end{align}
\end{subequations}
Eq. \eqref{p1:objective_funtion} corresponds to the objective function, which aims to minimize the cost of generation and battery discharge while minimizing the load shedding cost. Eq. in \eqref{p1:gen_limts} corresponds to generation limits, \eqref{p1:energy_balance} corresponds to the energy balance. Eq. in \eqref{p1:line_limits} corresponds to line flow limits, and eq. \eqref{p1:line_angles} to the angle-line relation. Eqs. \eqref{p1:battery_dynamics} and \eqref{p1:battery_limits} correspond to battery dynamics and limits, respectively. Eqs. \eqref{p1:battery_dis_limit} and \eqref{p1:battery_ch_limit} correspond to the battery discharge and charge limits. The variables on the left side correspond to the dual variables that will be discussed in Section \ref{subsection:dual}.

Note that, since ignition risk data are available at a daily granularity and representative weeks are considered, while the remaining system information is provided on an hourly basis, it is necessary to ensure that the battery dynamic constraint in \eqref{p1:battery_dynamics} maintains temporal consistency. This is achieved by applying the following logic:
\begin{subequations}
\begin{align}
    e_{s}^{t+1,d,w} = e_{s}^{t,d,w} + \eta c_{s}^{t+1,d,w} - \frac{1}{\eta} d_{s}^{t+1,d,w},\label{p1:battery_dynamics_hour}\\
    e_{s}^{1,d+1,w} = e_{s}^{T,d,w} + \eta c_{s}^{1,d+1,w} - \frac{1}{\eta} d_{s}^{1,d+1,w},\label{p1:battery_dynamics_days}\\
    e_{s}^{1,1,w+1} = e_{s}^{T,D,w} + \eta c_{s}^{1,1,w+1} - \frac{1}{\eta} d_{s}^{1,1,w+1}.\label{p1:battery_dynamics_week}
\end{align}
\end{subequations}

\subsection{Uncertainty set design} \label{section:uncertainty_sets}
An important aspect of the proposed robust optimization model is the design of the uncertainty sets $\mathcal{U}_{r}$ and $\mathcal{U}_{\alpha}$ for transmission line de-energization and renewable generation, respectively. In this paper, $\mathcal{U}_{r}$ described as a box uncertainty set as follows:
\begin{subequations}
\begin{align}
    \mathcal{U}_{r}(z^{\rm UG}) = & \Big\{z^{d,w}_{l},\gamma_{l}^{d,w}, r_{l}^{\star d,w}, r_{l}^{d,w}: \nonumber\\ 
    & \quad \gamma_{l}^{d,w} \in [0,1],~z^{d,w}_{l} \in \{0,1\},~ r^{\star d,w}_{l} \geq 0, \nonumber\\
    & r_l^{\text{avg},d,w} - \gamma_{l}^{d,w} r_l^{\text{dev},d,w} \leq r_{l}^{d,w}, \label{p2:risk_bounds_1}\\
    & r_{l}^{d,w} \leq r_l^{\text{avg},d,w} + \gamma_{l}^{d,w} r_l^{\text{dev},d,w}, \label{p2:risk_bounds_2}\\
    & \sum_{l} \gamma_{l}^{d,w} \leq \Gamma^{d,w} \sqrt{|L|}, \label{p2:budget}\\
    & r^{\star d,w}_{l} \leq M~z_{l}^{d,w} + r^{\rm th}_{l}, \label{p2:risk_threshold_up}\\
    & -M~(1-z_{l}^{d,w})+ r^{\rm th}_{l} \leq r^{\star d,w}, \label{p2:risk_threshold_dn}\\
    & r^{\star d,w}_{l} \leq r^{d,w}_{l}, \label{p2:real_risk}\\
    & r^{\star d,w}_{l} \leq (1-z^{\rm UG}_{l})~M \Big\} \label{p2:risk_undergroung},
\end{align}
\end{subequations}
where, $r_l^{\text{avg},d,w}$ and $r_l^{\text{dev},d,w}$ are the historical line ignition risk mean and deviation, respectively. The binary variable $z^{d,w}_{l}$ associated with the transmission line $l$ is equal to $1$ if the ignition risk exceeds its risk threshold $r^{\rm th}_{l}$ and must be de-energized and $0$ otherwise, $\gamma_{l}^{d,w}$ is the risk deviation associated with $l$ from its mean ignition risk value, $r_{l}^{d,w}$ is the risk perceived by the line $l$ due to its location, and $r_{l}^{\star d,w}$ is the real ignition risk of $l$ that takes into account if the line is underground. Thus, eqs. \eqref{p2:risk_bounds_1} and \eqref{p2:risk_bounds_2} correspond to line risk limits based on their historical information. Eq. \eqref{p2:budget} corresponds to the maximum risk deviation limit based on the budget parameter $\Gamma^{d,w}$. Eqs. \eqref{p2:risk_threshold_up} and \eqref{p2:risk_threshold_dn} correspond to the line de-energization logic based on the line's threshold risk. This is a simple assumption where each transmission line is de-energized if the real ignition risk exceeds their historical operational threshold. Finally, eqs. \eqref{p2:real_risk} and \eqref{p2:risk_undergroung} define the real ignition risk limits. Specifically, the limit is set to $0$ if the line is underground, and to the perceived ignition risk otherwise. Consequently, $\mathcal{U}(z^{\rm UG})$ depends on $z^{\rm UG}$, since the mitigation action of de-energization is only needed for lines that are not undergrounded.

For variable generation availability, we consider the uncertainty set described in \cite{verastegui2019adaptive, tapia2021robust}. In this paper, we consider both wind and solar generation as uncertain generation technologies, where their uncertainty set is described as follows.
\begin{subequations}
\begin{align}
    \mathcal{U}_{\alpha} = \Big\{\alpha_{r}^{t,d,w}:~& \alpha_{r}^{t,d,w} \in [\alpha_{r}^{\text{avg},t,d,w} \pm \alpha_{r}^{\text{dev},t,d,w}], \label{p7:risk_bounds}\\
    & \sum_{r} \frac{|\alpha_{r}^{t,d,w} -  \alpha_{r}^{\text{avg},t,d,w}|}{\alpha_{r}^{\text{dev},t,d,w}} \leq \Gamma^{t,d,w}_{\alpha} \sqrt{|W|}, \label{p7:budget} \Big\}
\end{align}
\end{subequations}
Where, eq. \eqref{p7:risk_bounds} defines the limits of the capacity factor, and eq. \eqref{p7:budget} corresponds to the maximum deviation limit of the capacity factor with respect to its mean value and bounded by the budget parameter $\Gamma^{t,d,w}_{\alpha}$.

\subsection{Solution approach for the worst-case wildfire ignition risk and variable generation availability} \label{subsection:dual}
In order to solve the $\max\min$ structure in \eqref{p0:invest_problem}, we reformulate the lower-level minimization problem using strong duality (and the dual variables in parentheses in \eqref{p1:economic_dispatch}), obtaining a pure maximization problem. As a result, we obtain the following equivalent bilinear optimization problem:
\begin{subequations} \label{p3:dual}
\begin{align}
    \max_{\substack{\alpha, z\\ \overline{\Delta}, \underline{\Delta}, \Delta, \pi, \mu}} & \Bigg( \sum_{t,d,w} \Big[ \sum_{l} (\underline{\Delta}^{t,d,w}_{l}+\overline{\Delta}^{t,d,w}_{l})~\big(-\overline{f}_{l} (1-z_{l}^{d,w})\big) \nonumber\\
    & + \sum_{l} M z^{d,w}_{l} (\overline{\Theta}_{l}^{t,d,w} + \underline{\Theta}_{l}^{t,d,w}) \nonumber \\
    & + \sum_{b} \pi^{t,d,w}_{b} d^{t,d,w}_{b} + \sum_{g} \big(\underline{\Delta}^{t,d,w}_{g} -\overline{\Delta}^{t,d,w}_{g} \big)p_g^{\rm max} \nonumber \\
    & + \sum_{r} -\overline{\Delta}^{t,d,w}_{r} \alpha_{r}^{t,d,w} p_{r}^{\rm max} - \sum_{s} (\overline{\Theta}^{t,d,w}_{s} + \overline{\Phi}^{t,d,w}_{s} )P_{s} \nonumber\\ 
    & + \sum_{s} (\underline{\Delta}^{t,d,w}_{s} X^{\min}_{s} - \overline{\Delta}^{t,d,w}_{s} X^{\max}_{s}) \nonumber\\
    & + \sum_{s} \mu^{1,d,w}_{s} E_{s}^{1,d,w} + \mu_{s}^{T,d,w} E_{s}^{T,d,w} \Big] \Bigg) \label{p3:objective_function}\\
    \text{s.t.}\quad& \overline{\Delta}^{t,d,w}_{l} - \underline{\Delta}^{t,d,w}_{l} + \overline{\Theta}^{t,d,w}_{l} - \underline{\Theta}^{t,d,w}_{l} \nonumber  \label{p3:init}\\
    & \qquad \qquad \qquad - \pi^{t,d,w}_{b^{\rm in}(l)} + \pi^{t,d,w}_{b^{\rm out}(l)} = 0,\\
    & C_{g} + \overline{\Delta}^{t,d,w}_{g} - \underline{\Delta}^{t,d,w}_{g} - \pi^{t,d,w}_{b(g)} \geq 0,\\
    & \overline{\Delta}^{t,d,w}_{r} - \pi^{t,d,w}_{b(r)} \geq 0,\\
    & \sum_{l|b = b^{\rm in}(l)} \overline{\Theta}^{t,d,w}_{l} B_l - \sum_{l|b = b^{\rm out}(l)} \underline{\Theta}^{t,d,w}_{l} B_l = 0,\\
    & C^{LS}_{b} - \pi^{t,d,w}_{b} \geq 0,\\
    & -C_{s} + \pi_{b(s)}^{t,d,w} - \frac{1}{\eta}\mu_{s}^{t,d,w} - \underline{\Theta}_{s}^{t,d,w} + \overline{\Theta}_{s}^{t,d,w} \geq 0,\\
    & \pi_{b(s)}^{t,d,w} + \eta \mu_{s}^{t,d,w} - \underline{\Phi}_{s}^{t,d,w} + \overline{\Phi}_{s}^{t,d,w} \geq 0,\\
    & -\mu_{s}^{t,d,w} + \mu_{s}^{t+1,d,w} + \underline{\Delta}_{s}^{t,d,w} -\overline{\Delta}_{s}^{t,d,w} \geq 0,\\
    & \overline{\Delta}^{t,d,w}_{l},~ \underline{\Delta}^{t,d,w}_{l},~ \overline{\Delta}^{t,d,w}_{g},~ \underline{\Delta}^{t,d,w}_{g},~ \overline{\Delta}^{t,d,w}_{r} \geq 0,\\
    & \overline{\Delta}^{t,d,w}_{s},~ \underline{\Delta}^{t,d,w}_{s},~ \overline{\Phi}^{t,d,w}_{s},~ \underline{\Phi}^{t,d,w}_{s} \geq 0,\\
    & \overline{\Theta}^{t,d,w}_{s},~ \underline{\Theta}^{t,d,w}_{s},~ \overline{\Theta}^{t,d,w}_{l},~ \underline{\Theta}^{t,d,w}_{l} \geq 0\\
    & \pi^{t,d,w}_{b}, \mu_{s}^{t,d,w} = \text{free},\\
    & \alpha_r^{t,d,w} \in \mathcal{U}_{\alpha}, \quad z_l^{d,w} \in \mathcal{U}_{r} \label{p3:end}
\end{align}
\end{subequations}
where, $\overline{\Delta}^{t,d,w}_{l}$ and $\underline{\Delta}^{t,d,w}_{l}$ are the dual variables associated with the transmission line limits in \eqref{p1:line_limits}, $\overline{\Theta}^{t,d,w}_{l}$ and $\underline{\Theta}^{t,d,w}_{l}$ are the dual variables associated with the angle phase in \eqref{p1:line_angles}. While $\overline{\Delta}^{t,d,w}_{g}$ and $\underline{\Delta}^{t,d,w}_{g}$ are the dual variables associated with the thermal generation limits in \eqref{p1:gen_limts}, $\overline{\Delta}^{t,d,w}_{r}$ is the dual variable associated with the maximum variable generation limit \eqref{p1:gen_renewable_limts}. Similarly, $\overline{\Delta}^{t,d,w}_{s}$ and $\underline{\Delta}^{t,d,w}_{s}$, $\overline{\Theta}^{t,d,w}_{s}$ and $\underline{\Theta}^{t,d,w}_{s}$, and $\overline{\Phi}^{t,d,w}_{s}$ and $\underline{\Phi}^{t,d,w}_{s}$ are the dual variables associated with the energy storage limits in \eqref{p1:battery_limits}, the limits of charge in \eqref{p1:battery_dis_limit} and the limits of discharge in \eqref{p1:battery_ch_limit}, respectively. $\pi^{t,d,w}_{b}$ and $\mu_{s}^{t,d,w}$ correspond to the dual variables associated with \eqref{p1:energy_balance} and \eqref{p1:battery_dynamics}, respectively. Finally, finally, $\alpha_{r}^{t,d,w}$ and $z_{l}^{d,w}$ are the variable generation availability and transmission line de-energization binary variable from Section \ref{section:uncertainty_sets}.

The objective function in eq. \eqref{p3:objective_function} contains a bilinear term, which corresponds to the multiplication between the availability of the renewable resource $\alpha_{r}^{t,d,w}$, the dual variable $\overline{\Delta}^{t,d,w}_{r}$. In addition, eq. \eqref{p3:objective_function} contains the multiplication between the continuous variables $\underline{\Delta}^{t,d,w}_{l}$ and $\overline{\Delta}^{t,d,w}_{l}$ with the integer variable $z_{l}^{d,w}$. We implement the traditional big-M method to reformulate the multiplication between a continuous and an integer variable as follows:
\begin{subequations}
\begin{align} 
    & \underline{\nu}_{l}^{t,d,w} \geq -M(1-z_{l}^{d,w}) + \underline{\Delta}_{l}^{t,d,w}, \label{p7:init}\\
    & \underline{\nu}_{l}^{t,d,w} \leq M(1-z_{l}^{d,w}) + \underline{\Delta}_{l}^{t,d,w},\\
    & -Mz_{l}^{d,w} \leq \underline{\nu}_{l}^{t,d,w} \leq Mz_{l}^{d,w},\\
    & \overline{\nu}_{l}^{t,d,w} \geq -M(1-z_{l}^{d,w}) + \overline{\Delta}_{l}^{t,d,w},\\
    & \overline{\nu}_{l}^{t,d,w} \leq M(1-z_{l}^{d,w}) + \overline{\Delta}_{l}^{t,d,w},\\
    & -Mz_{l}^{d,w} \leq \overline{\nu}_{l}^{t,d,w} \leq Mz_{l}^{d,w},\\
    & \overline{\nu}_{l}^{t,d,w}, \underline{\nu}_{l}^{t,d,w} \geq 0. \label{p7:end}
\end{align}
\end{subequations}
Likewise, for the multiplication between continuous variables $\underline{\Theta}^{t,d,w}_{l}$ and $\overline{\Theta}^{t,d,w}_{l}$ with the integer variable $z_{l}^{d,w}$, we reformulate them as follows:
\begin{subequations}
\begin{align}
    & \underline{\zeta}_{l}^{t,d,w} \geq -M(1-z_{l}^{d,w}) + \underline{\Theta}_{l}^{t,d,w}, \label{p8:init}\\
    & \underline{\zeta}_{l}^{t,d,w} \leq M (1-z_{l}^{d,w}) + \underline{\Theta}_{l}^{t,d,w},\\
    & -Mz_{l}^{d,w} \leq \underline{\zeta}_{l}^{t,d,w} \leq Mz_{l}^{d,w},\\
    & \overline{\zeta}_{l}^{t,d,w} \geq -M (1-z_{l}^{d,w}) + \overline{\Theta}_{l}^{t,d,w},\\
    & \overline{\zeta}_{l}^{t,d,w} \leq M (1-z_{l}^{d,w}) + \overline{\Theta}_{l}^{t,d,w},\\
    & -Mz_{l}^{d,w} \leq \overline{\zeta}_{l}^{t,d,w} \leq Mz_{l}^{d,w},\\
    & \overline{\zeta}_{l}^{t,d,w}, \underline{\zeta}_{l}^{t,d,w} \geq 0. \label{p8:end}
\end{align}
\end{subequations}

In the case of the bilinear term, we implement the binary expansion method in \cite{gupte2013solving} to discretize  $\alpha_{r}^{t,d,w}$ and use the big-M method to obtain an approximating mixed-integer optimization problem that can be solved using off-the-shelf solvers \cite{tapia2021robust}. Taking minimum battery capacity size as $X_{s}^{\rm min} = 0$, the resulting bilinear terms are shown below:
\begin{subequations}
\begin{align}
    \alpha_{r}^{t,d,w} & = \sum_{n} \frac{1}{2^{n-1}} u_{r,n}^{t,d,w} \alpha_{r}^{\text{dev},t,d,w} \nonumber\\
    & + (\alpha_{r}^{\text{avg},t,d,w} - \alpha_{r}^{\text{dev},t,d,w}), \label{eq:binary_eq}\\
    \overline{\Delta}_{r}^{t,d,w} \alpha_{r}^{t,d,w} & = \overline{\Delta}_{r}^{t,d,w} \sum_{n} \frac{1}{2^{n-1}} u_{r,n}^{t,d,w} \alpha_{r}^{\text{dev},t,d,w} \nonumber \\
    & + \overline{\Delta}_{r}^{t,d,w} (\alpha_{r}^{\text{avg},t,d,w} - \alpha_{r}^{\text{dev},t,d,w}),\\
    & = \alpha_{r}^{\text{dev},t,d,w} \sum_{n} \frac{1}{2^{n-1}}  \overline{\varphi}_{s,n}^{t,d,w} \nonumber\\
    & + \overline{\Delta}_{r}^{t,d,w} (\alpha_{r}^{\text{avg},t,d,w} - \alpha_{r}^{\text{dev},t,d,w}),
\end{align}
\end{subequations}
where $\overline{\varphi}$ is constrained as follows:
\begin{subequations}
\begin{align}
    & -M(1-u_{r,n}^{t,d,w}) + \overline{\Delta}_{r}^{t,d,w} \leq \overline{\varphi}_{r,n}^{t,d,w}, \label{p10:init}\\
    & \overline{\varphi}_{r,n}^{t,d,w} \leq M(1-u_{r,n}^{t,d,w}) + \overline{\Delta}_{r}^{t,d,w},\\
    & -Mu_{r,n}^{t,d,w} \leq \overline{\varphi}_{r,n}^{t,d,w} \leq Mu_{r,n}^{t,d,w},\\
    & \overline{\varphi}_{r,n}^{t,d,w} \geq 0, \label{p10:end}
\end{align}
\end{subequations}
where, $n$ represents the number of segments of the binary expansion. Parameters $ \alpha_{r}^{\text{avg},t,d,w}$ and $ \alpha_{r}^{\text{dev},t,d,w}$ are the mean and deviation of the renewable availability factor, respectively. Increasing the value of $n$ improves the \textit{resolution} of $\alpha^{t,d,w}_{r}$, but increases the number of binary variables $u_{r,n}^{t,d,w}$ and the computational burden for the proposed model. Then, repacing eq. \eqref{eq:binary_eq} in eq. \eqref{p7:budget}, we approximate the budget constraint for the binary expansion of the capacity factor as follows:
\begin{equation}
    \sum_{r} \Big|1 - \sum_{n} \frac{1}{2^{n-1}} u_{r,n}^{t,d,w}\Big| \leq \Gamma^{t,d,w}_{\alpha} \sqrt{|W|} \label{p11}
\end{equation}

\noindent Thus, the $\max\min$ model could be expressed as follows:
\begin{subequations} \label{p3:dual2}
\begin{align}
    \max_{\alpha, z}~& \max_{\overline{\Delta}, \underline{\Delta}, \Delta, \pi, \mu} \tag{12} \\
    & \Bigg( \sum_{t,d,w} \Big[ \sum_{l} (\underline{\Delta}^{t,d,w}_{l}+\overline{\Delta}^{t,d,w}_{l})~\big(-\overline{f}_{l} \big) + (\underline{\nu}_{l}^{t,d,w} + \overline{\nu}_{l}^{t,d,w})\nonumber\\
    & + \sum_{l} M (\overline{\zeta}_{l}^{t,d,w} + \underline{\zeta}_{l}^{t,d,w}) - \sum_{s} (\overline{\Theta}^{t,d,w}_{s} + \overline{\Phi}^{t,d,w}_{s} )P_{s} \nonumber\\ 
    & + \sum_{s} (\underline{\Delta}^{t,d,w}_{s} X^{\min}_{s} - \overline{\Delta}^{t,d,w}_{s} X^{\max}_{s}) \nonumber\\
    & + \sum_{b} \pi^{t,d,w}_{b} d^{t,d,w}_{b} + \sum_{g} \big(\underline{\Delta}^{t,d,w}_{g} -\overline{\Delta}^{t,d,w}_{g} \big)p_g^{\rm max} \nonumber \\
    & - \sum_{r} p_{r}^{\rm max} \Big( \alpha_{r}^{\text{dev},t,d,w} \sum_{n} \frac{1}{2^{n-1}}  \overline{\varphi}_{s,n}^{t,d,w} \nonumber\\
    & + \overline{\Delta}_{r}^{t,d,w} (\alpha_{r}^{\text{avg},t,d,w} - \alpha_{r}^{\text{dev},t,d,w}) \Big) \nonumber\\
    & + \sum_{s} \mu^{1,d,w}_{s} E_{s}^{1,d,w} + \mu_{s}^{T,d,w} E_{s}^{T,d,w} \Big] \Bigg) \nonumber\\
    \text{s.t.}\quad& \eqref{p3:init}-\eqref{p3:end},~\eqref{p7:init}-\eqref{p7:end},~\eqref{p8:init}-\eqref{p8:end} \nonumber\\
    & \eqref{p10:init}-\eqref{p10:end},~\eqref{p11}. \nonumber
\end{align}
\end{subequations}

\section{Solution approach for ignition-aware robust CEP model}
In order to solve the robust optimization model proposed in \eqref{p0:invest_objective}, we conveniently reformulate it as follows.
\begin{subequations} \label{p4:invest}
\begin{align} 
    \min_{X,z^{\rm UG}} ~& \sum_{l} C_{l}^{\rm UG} z_{l}^{UG} + \sum_{s} C_{s}^{X}~X_{s}^{\rm max} + \delta, \label{p4:invest_objective}\\
    \text{s.t.}\quad &0 \leq X_{s}^{\rm max} \leq X_{s}^{\rm sup}, \label{p4:investment_limits}\\
    & \delta \geq \min_{y \in Y(X,z^{\rm UG},\xi)} f(\cdot) \quad \forall~\xi \in \Xi, \label{p4:wost-case}\\
    & z_{l}^{\rm UG} \in \{0,1\},    
\end{align}
\end{subequations}
where, $\xi$ represents the realization of uncertainty within the set $\Xi = \mathcal{U}_r(z^{\rm UG})\times\mathcal{U}_{\alpha}$, while $y = (p,f,e,d,c)$ and $Y$ is the feasible set for the dispatch variables under the realization of $\xi$ and the investment $X,z^{\rm UG}$. Finally, $\delta$ corresponds to the worst-case operational cost in the objective function in \eqref{p3:dual2} represented as $f(\cdot)$. Thus, eq. \eqref{p4:investment_limits} represents the maximum investment in storage capacity for candidate projects and eq. \eqref{p4:wost-case} ensures that $\delta$ captures the dispatch cost of the worst-case both uncertainties realization, ignition risk and variable generation availability. We use a column constraint and generation method in \cite{verastegui2019adaptive} to solve the ignition-aware robust CEP problem.

\section{Case Study}
In this section, we evaluate the proposed robust CEP model using a representation of the San Diego power system. We analyze expansion plans under varying levels of uncertainty in wildfire ignition risk and renewable generation availability. Two investment schemes are considered: 
\begin{itemize}
    \item [(i)] \textbf{Scheme 1}: only storage investment is allowed;
    \item [(ii)] \textbf{Scheme 2}: investments in storage and underground transmission lines are allowed.
\end{itemize} 

\noindent The cost of load shedding is set at $C^{\rm LS} = \$20{,}000/\text{MWh}$. The cost of undergrounding a transmission line is $C^{\rm UG} = \$7{,}000{,}000/\text{mile}$ \cite{piansky2025optimizing} and the line lengths are estimated from the coordinates of the associated nodes. The cost of battery capacity is $C^X = \$1{,}000{,}000/\text{MW}$ of installed capacity, assuming a battery efficiency of \( \eta = 0.95 \). A maximum battery size of 400~MW per node is considered. All installation costs are annualized based on the expected lifetime of each investment asset.

We use the maximum ignition risk method from \cite{piansky2024quantifying} to assign wildfire risk levels to each transmission line. The mean and maximum deviations of ignition risk are estimated using the approach described in Section \ref{section:ignition_maps}, based on 2024 data from \cite{usgs_wfpi_2025} and a selected subset of three representative weeks per year. The corresponding representative-week weights are $[0.583, 0.400, 0.017]$, indicating that the first two weeks are dominant, accounting for 98.3\% of the total weight, while the third week represents a low-frequency condition with only a marginal effect on the simulations. De-energization thresholds for each line are determined using the 95th percentile of its historical ignition risk exposure.

We employ a test case consisting of a 24-bus representation of the San Diego power system, which includes 38 transmission lines, 33 existing conventional generation units, 3 existing variable generation units (that consider wind and solar generation), and 8 energy storage systems (3 existing and 5 candidate units). Demand profiles are based on 365 days of the base year 2026. This system is based on the IEEE 24-bus test system for transmission lines and generators data, and is mapped to the San Diego region in CATS \cite{taylor2023california} to incorporate its georeferenced information, demand profiles.

All numerical results presented below are on an annual basis and the maximum number of iterations is set to 20. The model is implemented in \textit{Pyomo} and solved using the \textit{Gurobi} optimizer.

\subsection{Results for the San Diego case study}
We run a full-year case using three representative weeks. The resulting investment decisions, battery sizes and the number of underground lines are summarized in Table \ref{Table:Scheme1} and \ref{Table:Scheme2} for Scheme 1 and 2, respectively.

\vspace{-1.5mm}
\begin{table}[htbp]
\centering
\caption{Battery investment under investment Scheme 1}
\vspace{-1.5mm}
\begin{tabular}{l|c|c|c}
\toprule
& \multicolumn{3}{c}{\textbf{Battery [MW]}}\\
$\Gamma_r$ & \multicolumn{1}{c|}{$\Gamma_{\alpha} = 0.1$} & \multicolumn{1}{c|}{$\Gamma_{\alpha} = 0.5$}& \multicolumn{1}{c}{$\Gamma_{\alpha} = 1.0$}\\
\midrule
\textbf{$0.1$} & 226.1 & 245.6 & 273.3\\
\textbf{$0.5$} & 252.5 & 194.1 & 183.9\\
\textbf{$1.0$}  & 209.4 & 272.8  & 229.3\\
\bottomrule
\end{tabular}
\vspace{-1.5mm}
\label{Table:Scheme1}
\end{table}

\begin{table}[htbp]
\vspace{-1.5mm}
\centering
\caption{Battery and underground (UG) line investment under investment Scheme 2}
\vspace{-1.5mm}
\begin{tabular}{l|cc|cc|cc}
\toprule
& \textbf{Battery} & \textbf{\# Lines} & \textbf{Battery} & \textbf{\# Lines} & \textbf{Battery} & \textbf{\# Lines}\\
& \textbf{[MW]} & \textbf{UG} & \textbf{[MW]} & \textbf{UG} & \textbf{[MW]} & \textbf{UG}\\
\textbf{$\Gamma_r$} & \multicolumn{2}{c|}{$\Gamma_{\alpha} = 0.1$} & \multicolumn{2}{c|}{$\Gamma_{\alpha} = 0.5$}& \multicolumn{2}{c}{$\Gamma_{\alpha} = 1.0$}\\
\midrule
\textbf{$0.1$} & 118.1 & 2 & 137.8 & 5 & 143.4 & 5\\
\textbf{$0.5$} & 137.1 & 2 & 145.6 & 4 & 155.3 & 6\\
\textbf{$1.0$}  & 157.2 & 3 & 172.7 & 3 & 188.0 & 6\\
\bottomrule
\end{tabular}
\vspace{-1.5mm}
\label{Table:Scheme2}
\end{table}

\begin{figure}[htbp!]
    \centering
    \begin{subcaptionbox}{Scheme 1: battery investment in [MW] for $\Gamma_\alpha = 0.1, ~\Gamma_r = 0.1 $\label{fig:subfig111}}[.32\linewidth]{\includegraphics[width=0.98\linewidth, trim={18mm, 4.5mm, 9mm, 8mm}, clip]{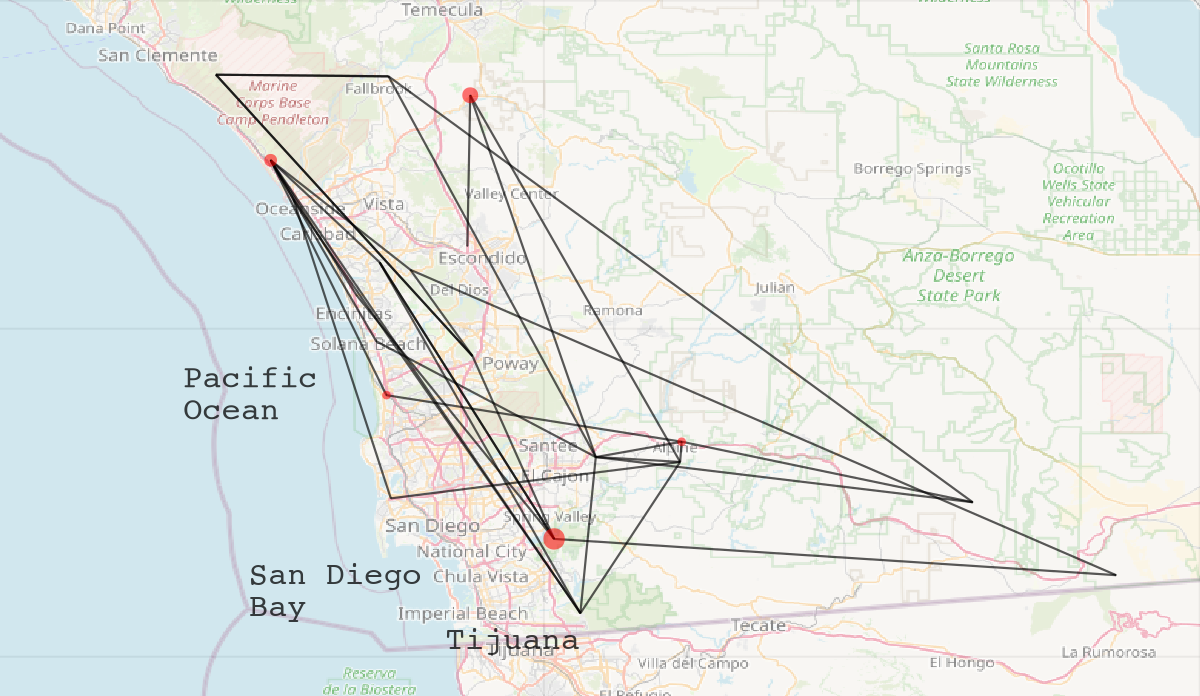}}
    \end{subcaptionbox}
    \begin{subcaptionbox}{Scheme 1: battery investment in [MW] for $\Gamma_\alpha = 0.1, ~\Gamma_r = 0.5 $\label{fig:subfig113}}[.32\linewidth]{\includegraphics[width=0.98\linewidth, trim={18mm, 4.5mm, 9mm, 8mm}, clip]{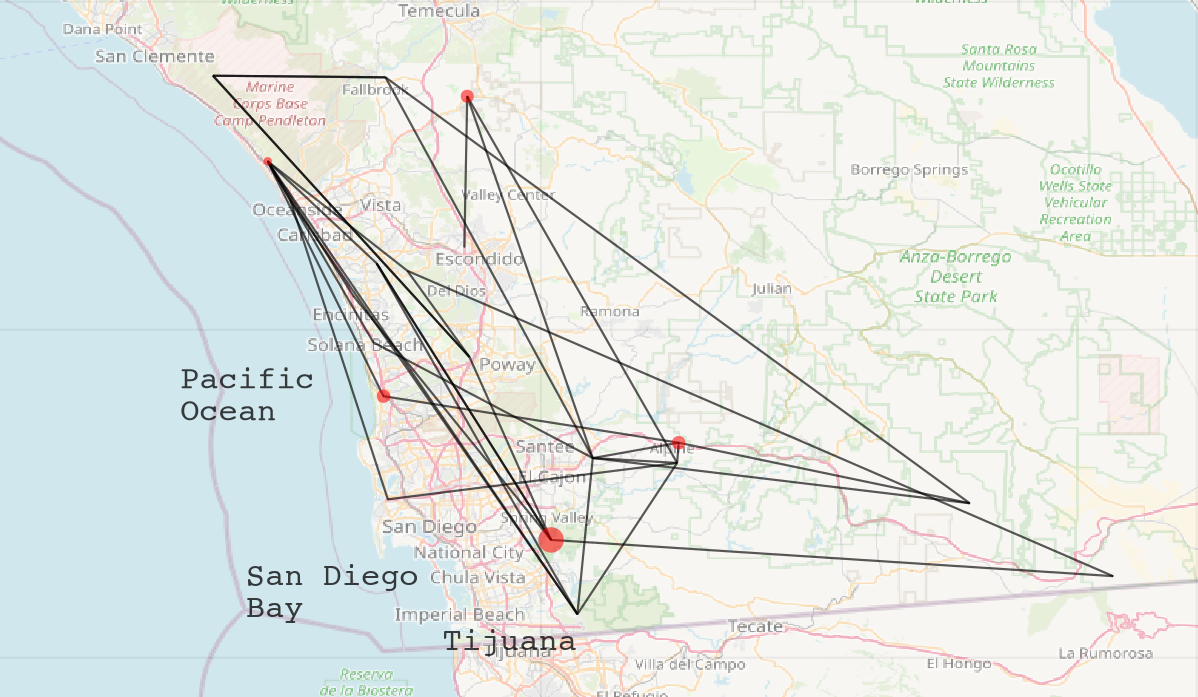}}
    \end{subcaptionbox}
    \begin{subcaptionbox}{Scheme 1: battery investment in [MW] for $\Gamma_\alpha = 0.1, ~\Gamma_r = 1.0 $\label{fig:subfig115}}[.32\linewidth]{\includegraphics[width=0.98\linewidth, trim={18mm, 4.5mm, 9mm, 8mm}, clip]{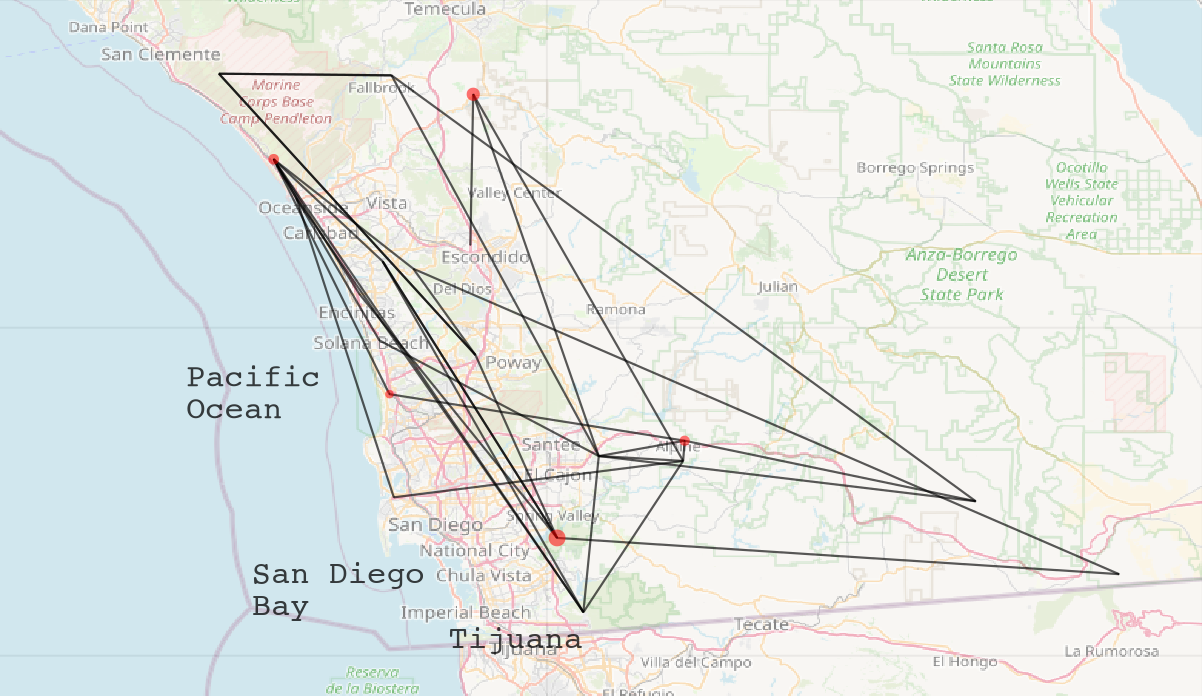}}
    \end{subcaptionbox}
    \hfill
    \begin{subcaptionbox}{Scheme 1: battery investment in [MW] for $\Gamma_\alpha = 1.0, ~\Gamma_r = 0.1 $\label{fig:subfig112}}[.32\linewidth]{\includegraphics[width=0.98\linewidth, trim={18mm, 4.5mm, 9mm, 8mm}, clip]{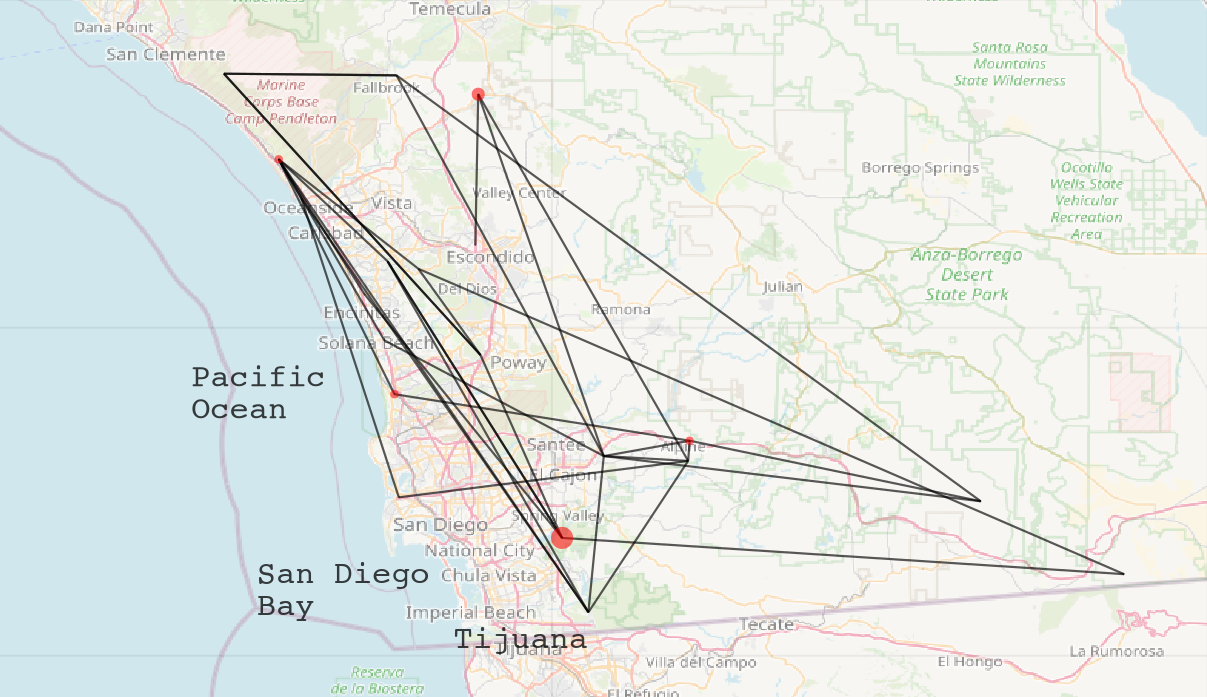}}
    \end{subcaptionbox}
    \begin{subcaptionbox}{Scheme 1: battery investment in [MW] for $\Gamma_\alpha = 1.0, ~\Gamma_r = 0.5 $\label{fig:subfig114}}[.32\linewidth]{\includegraphics[width=0.98\linewidth, trim={18mm, 4.5mm, 9mm, 8mm}, clip]{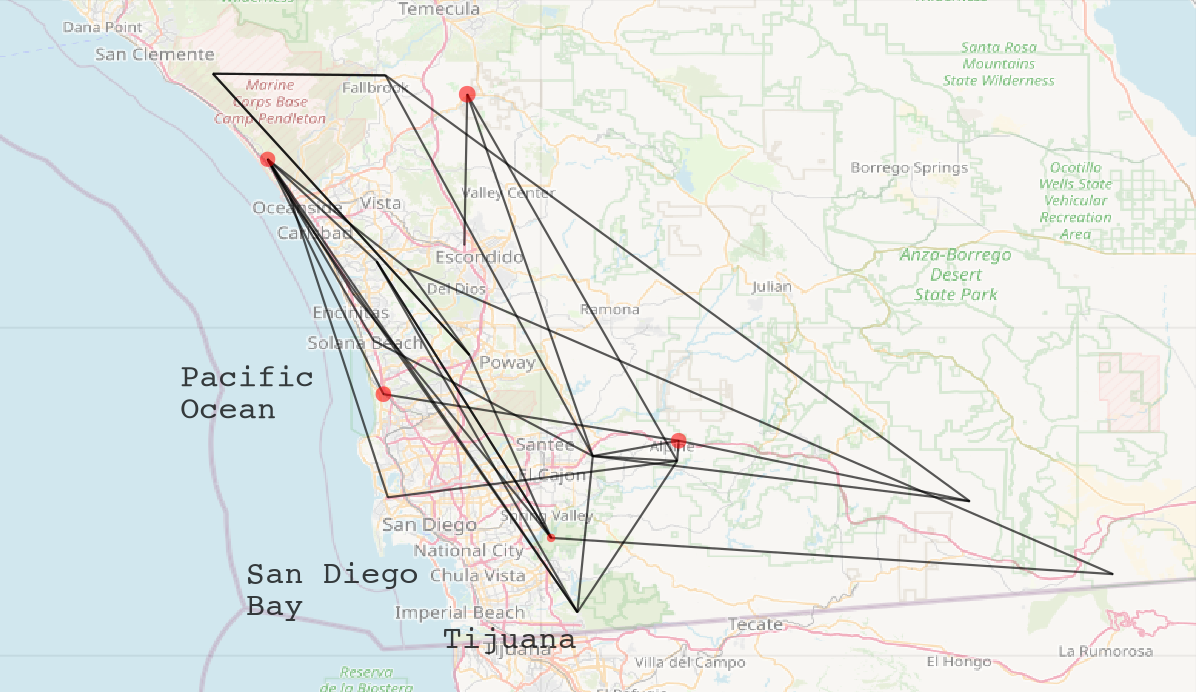}}
    \end{subcaptionbox}
    \begin{subcaptionbox}{Scheme 1: battery investment in [MW] for $\Gamma_\alpha = 1.0, ~\Gamma_r = 1.0 $\label{fig:subfig116}}[.32\linewidth]{\includegraphics[width=0.98\linewidth, trim={18mm, 4.5mm, 9mm, 8mm}, clip]{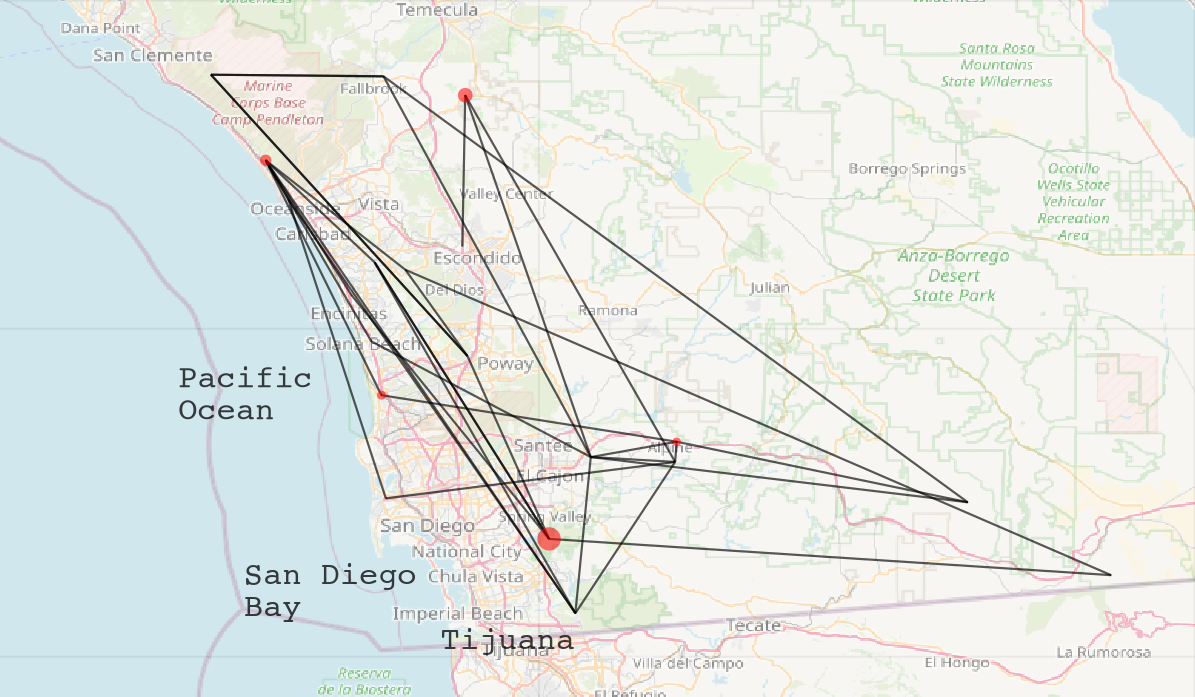}}
    \end{subcaptionbox}
    \vspace{-2.2mm}
    \caption{Storage capacity investment (red dots) solutions under investment Scheme 1. Largest circle equating to 100 MWh of installed capacity. The maps present a sensitivity on renewable $(\Gamma_\alpha)$ and line risk ($\Gamma_r)$ uncertainty sets as follow: (Top) $\Gamma_{\alpha}=0.1$; and (Bottom) $\Gamma_{\alpha}=1$;
    (Left) $\Gamma_r=0.1$; (Middle) $\Gamma_r=0.5$; (Right) $\Gamma_r=1$.}
    \label{fig:maps_scheme_1}
\end{figure}

\begin{figure}[h]
    \centering
    \begin{subcaptionbox}{Scheme 2: battery investment and lines UG for $\Gamma_\alpha = 0.1, ~\Gamma_r = 0.1 $\label{fig:subfig101}}[.32\linewidth]{\includegraphics[width=0.98\linewidth, trim={18mm, 4.5mm, 9mm, 8mm}, clip]{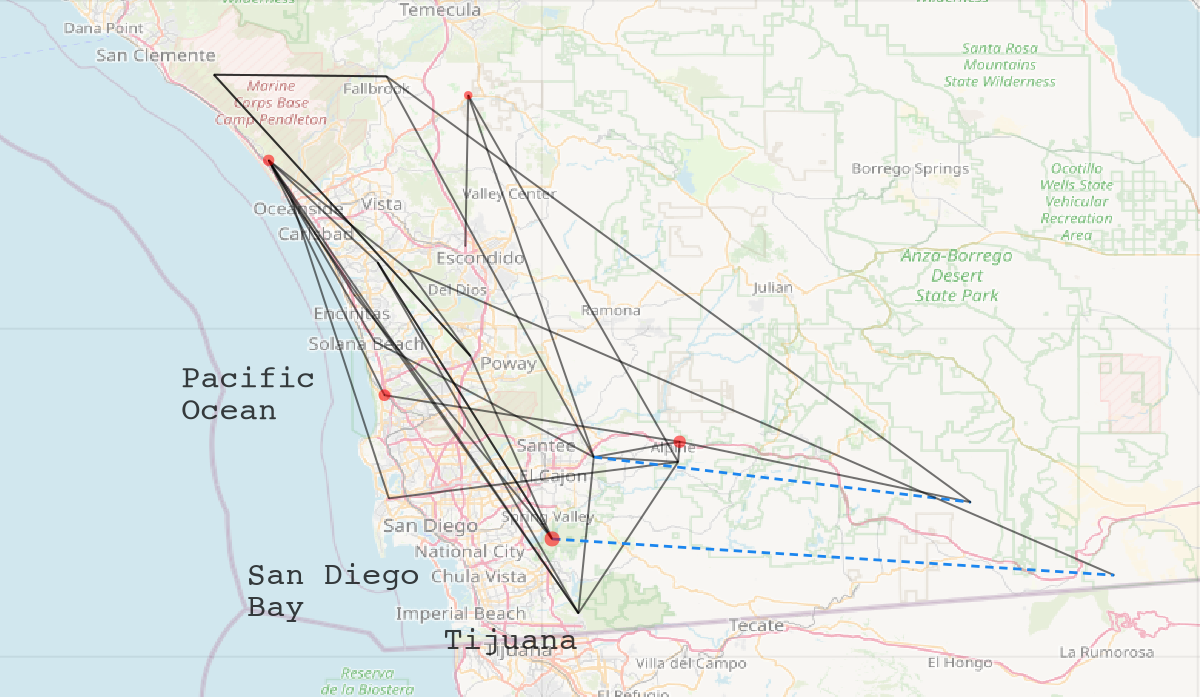}}
    \end{subcaptionbox}
    \begin{subcaptionbox}{Scheme 2: battery investment and lines UG for $\Gamma_\alpha = 0.1, ~\Gamma_r = 0.5 $\label{fig:subfig103}}[.32\linewidth]{\includegraphics[width=0.98\linewidth, trim={18mm, 4.5mm, 9mm, 8mm}, clip]{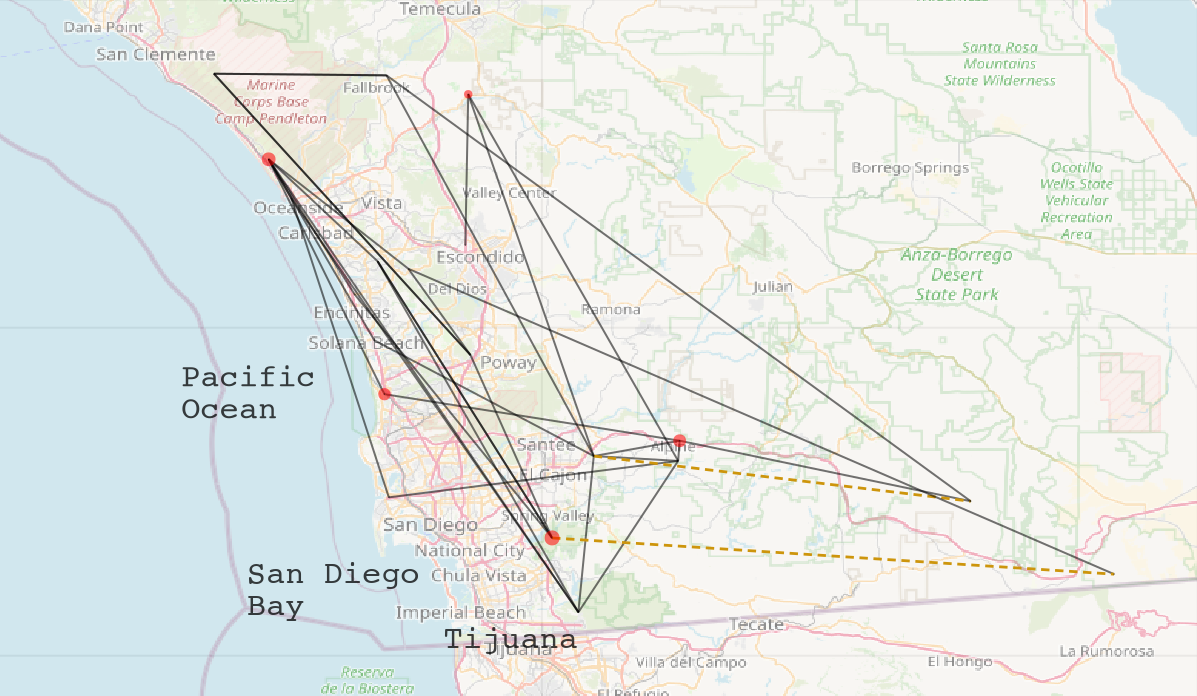}}
    \end{subcaptionbox}
    \begin{subcaptionbox}{Scheme 2: battery investment and lines UG for $\Gamma_\alpha = 0.1, ~\Gamma_r = 1.0 $\label{fig:subfig105}}[.32\linewidth]{\includegraphics[width=0.98\linewidth, trim={18mm, 4.5mm, 9mm, 8mm}, clip]{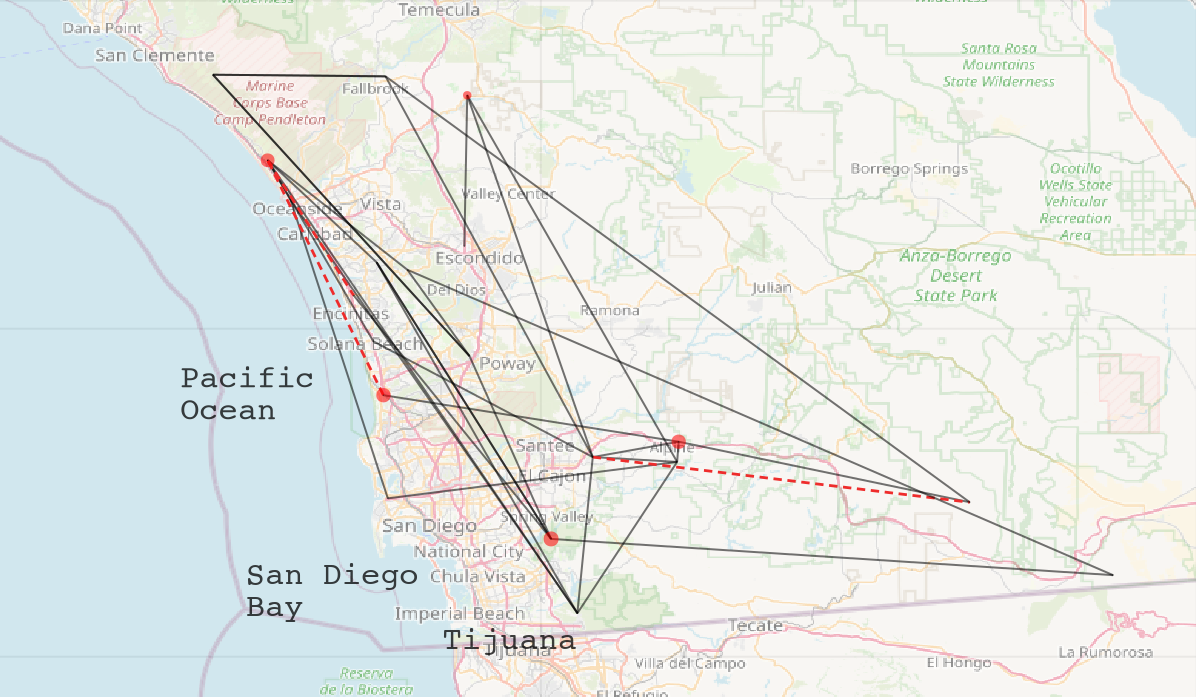}}
    \end{subcaptionbox}
    \hfill
    \begin{subcaptionbox}{Scheme 2: battery investment and lines UG for $\Gamma_\alpha = 1.0, ~\Gamma_r = 0.1 $\label{fig:subfig102}}[.32\linewidth]{\includegraphics[width=0.98\linewidth, trim={18mm, 4.5mm, 9mm, 8mm}, clip]{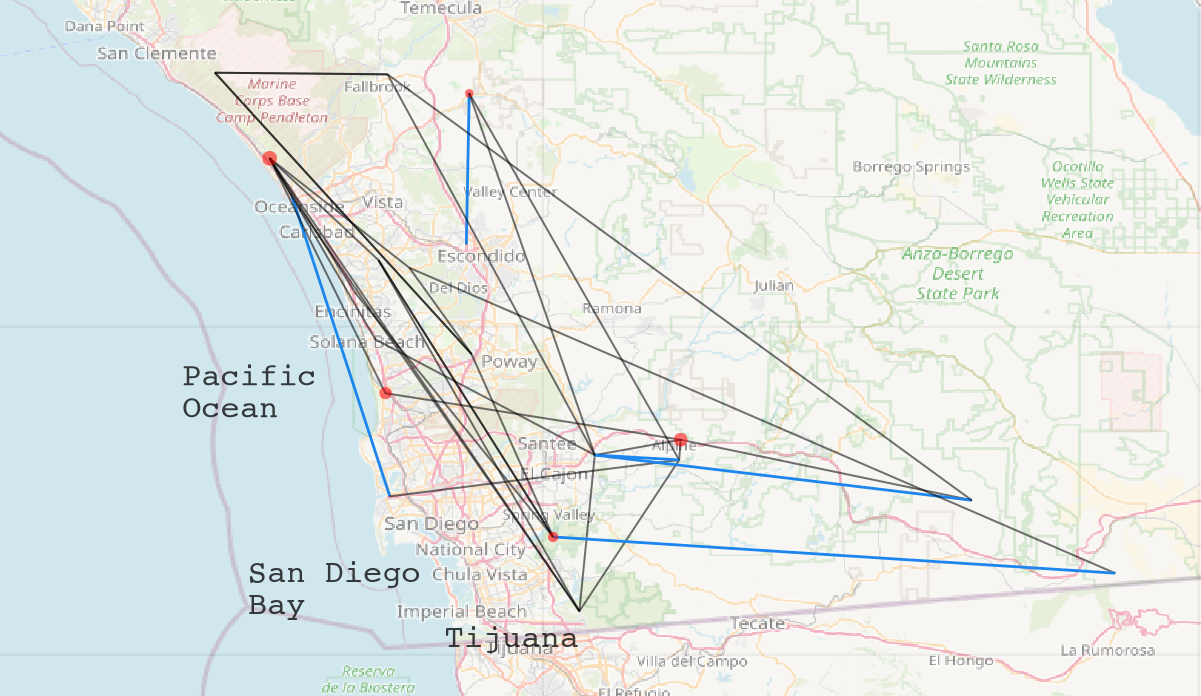}}
    \end{subcaptionbox}
    \begin{subcaptionbox}{Scheme 2: battery investment and lines UG for $\Gamma_\alpha = 1.0, ~\Gamma_r = 0.5 $\label{fig:subfig104}}[.32\linewidth]{\includegraphics[width=0.98\linewidth, trim={18mm, 4.5mm, 9mm, 8mm}, clip]{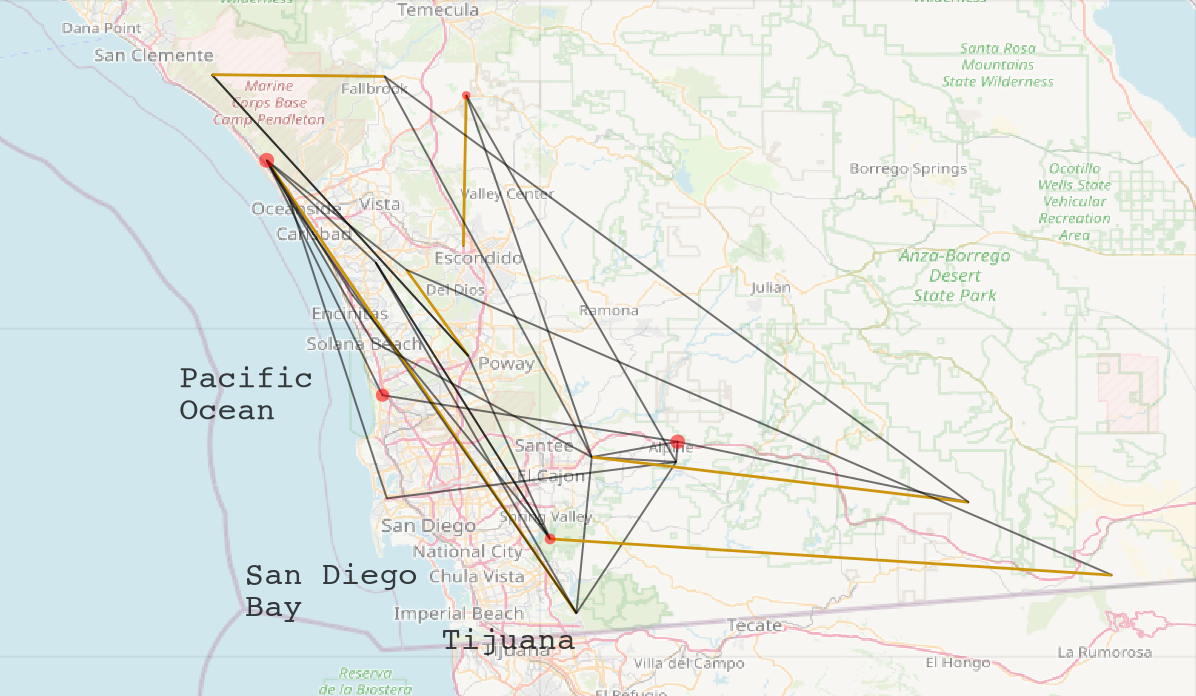}}
    \end{subcaptionbox}
    \begin{subcaptionbox}{Scheme 2: battery investment and lines UG for $\Gamma_\alpha = 1.0, ~\Gamma_r = 1.0 $\label{fig:subfig106}}[.32\linewidth]{\includegraphics[width=0.98\linewidth, trim={18mm, 4.5mm, 9mm, 8mm}, clip]{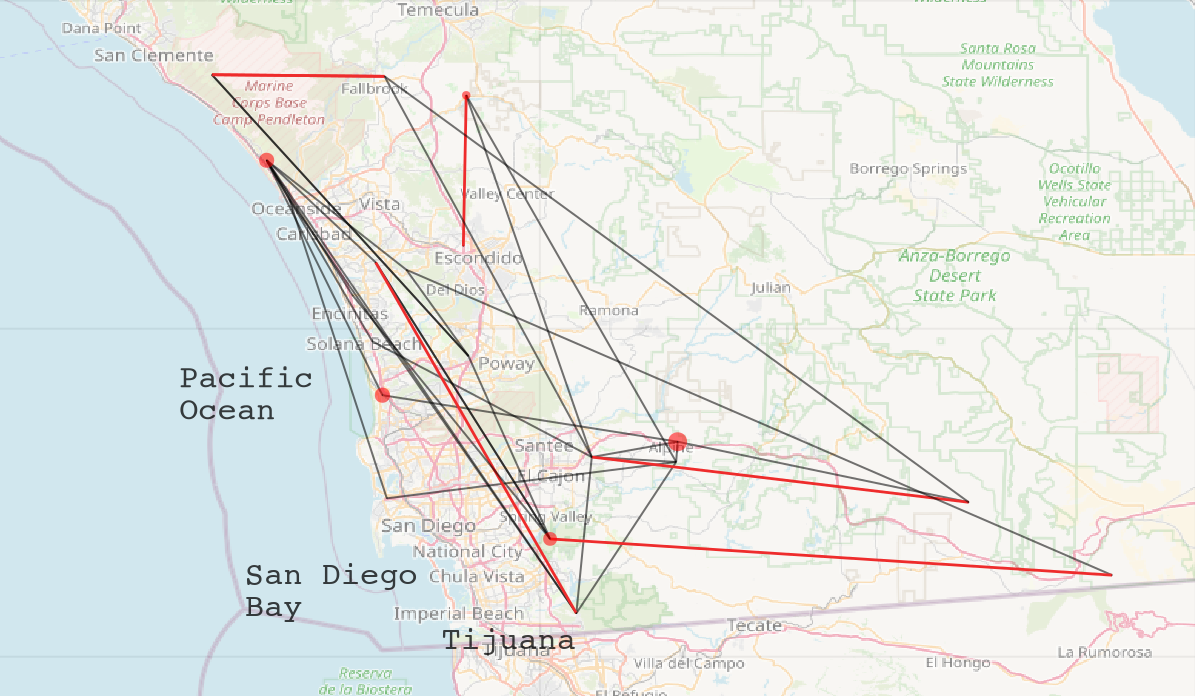}}
    \end{subcaptionbox}
    \vspace{-2mm}
    \caption{Line hardening (colored lines) and storage capacity investment (red dots) solutions under investment Scheme 2. Largest circle equating to 100 MWh of installed capacity. The maps present a sensitivity on renewable $(\Gamma_\alpha)$ and line wildfire ignition risk ($\Gamma_r)$ uncertainty sets as follow: (Top) $\Gamma_{\alpha}=0.1$; and (Bottom) $\Gamma_{\alpha}=1$; (Left) $\Gamma_r=0.1$; (Middle) $\Gamma_r=0.5$; (Right) $\Gamma_r=1$.}
    \label{fig:maps_scheme_2}
\end{figure}

The locational results are shown in Figure \ref{fig:maps_scheme_1} for Scheme 1 and Figure \ref{fig:maps_scheme_2} for Scheme 2. These maps illustrate the optimal placement batteries and undergrounded transmission lines across the San Diego power system. Red circles are proportional to the battery capacity installed at each bus, with the largest circle representing 150 MWh. Colored lines indicate transmission lines that have been undergrounded (hardened) for each of the sensitivity cases, while the remaining lines are shown in black. 

We observe that as $\Gamma_r$ increases, i.e. solution becomes more conservative, the number of undergrounded lines also increases and, the line selection changes due to the different impact of the different impact of line de-energization. In both cases, battery storage capacity tends to grow proportionally to increased conservatism of the solution. As $\Gamma_\alpha$ increases, the total battery storage investment increases too, but the distribution across nodes may shift. For example, the storage capacity at the southern node decreases from $41.6$ MWh to $23.1$ MWh when $\Gamma_\alpha$ changes from $0.1$ to $1.0$.

Table \ref{Table:load_shedding} compares the load shedding results between the two investment schemes. A greater reduction in load shedding is observed when both battery and underground line investments are considered. In all cases, daily load shedding is reduced from approx. $2.000$ MWh to $0$ MWh. It is important to note that the proposed model worst-case wildfire ignition risk allows for de-energization of any transmission lines, while in Scheme 1, the response is limited resulting in battery installations at only five nodes.

\begin{table}[htbp!]
\centering
\caption{Average load shedding per day [MWh]} 
\vspace{-1.5mm}
\begin{tabular}{l|ccc|ccc}
\toprule
& \multicolumn{3}{c|}{\textbf{Scheme 1}} & \multicolumn{3}{c}{\textbf{Scheme 2}}\\
$\Gamma_r$ & \multicolumn{1}{c|}{$\Gamma_{\alpha}\!=\!0.1$} & \multicolumn{1}{c|}{$\Gamma_{\alpha}\!=\!0.5$}& \multicolumn{1}{c|}{$\Gamma_{\alpha}\!=\!1$} & \multicolumn{1}{c|}{$\Gamma_{\alpha}\!=\!0.1$} & \multicolumn{1}{c|}{$\Gamma_{\alpha}\!=\!0.5$}& \multicolumn{1}{c}{$\Gamma_{\alpha}\!=\!1$}\\
\midrule
\textbf{$0.1$} & 2028 & 2002 & 2038 & 0 & 0 & 0\\
\textbf{$0.5$} & 2248 & 2259 & 2265 & 0 & 0 & 0\\
\textbf{$1.0$} & 2258 & 2924 & 2980 & 0 & 0 & 0\\
\bottomrule
\end{tabular}
\vspace{-1.5mm}
\label{Table:load_shedding}
\end{table}

\begin{figure}[htbp]
    \centering
    \begin{subcaptionbox}{$\Gamma_{\alpha}=1$, $\Gamma_r=1$.\label{fig:subfig71}}[.48\linewidth]
        {\includegraphics[width=\linewidth]{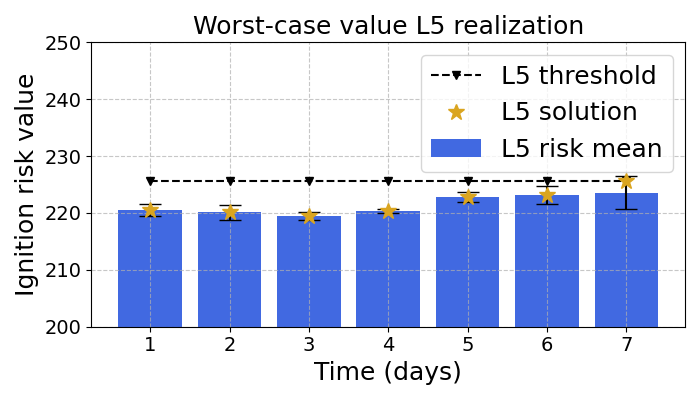}}
    \end{subcaptionbox}
    \hfill
    \begin{subcaptionbox}{$\Gamma_{\alpha}=1$, $\Gamma_r=0.1$.\label{fig:subfig72}}[.48\linewidth]
        {\includegraphics[width=\linewidth]{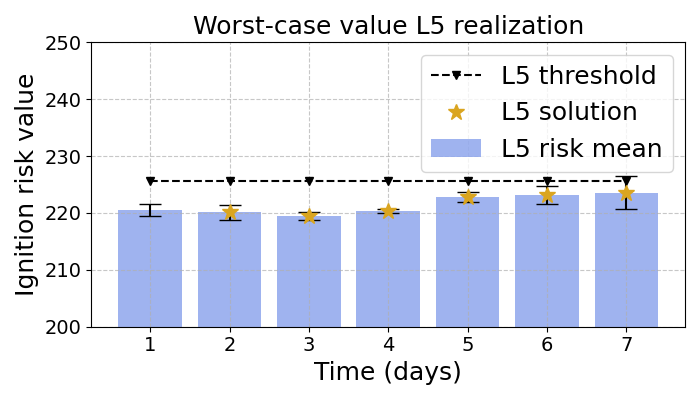}}
    \end{subcaptionbox}
    \begin{subcaptionbox}{$\Gamma_{\alpha}=1$, $\Gamma_r=1$.\label{fig:subfig73}}[.48\linewidth]
        {\includegraphics[width=\linewidth]{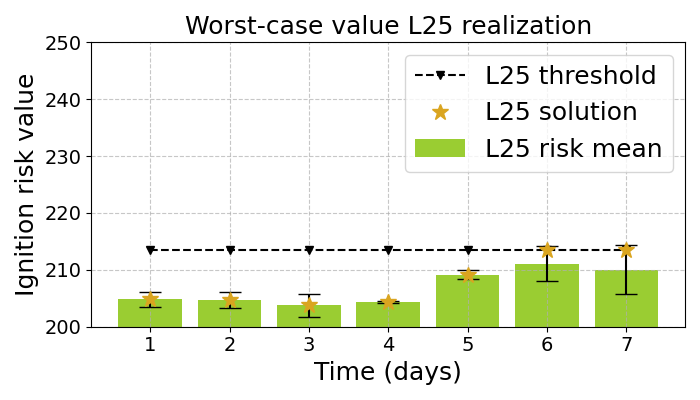}}
    \end{subcaptionbox}
    \hfill
    \begin{subcaptionbox}{$\Gamma_{\alpha}=1$, $\Gamma_r=0.1$. \label{fig:subfig74}}[.48\linewidth]
        {\includegraphics[width=\linewidth]{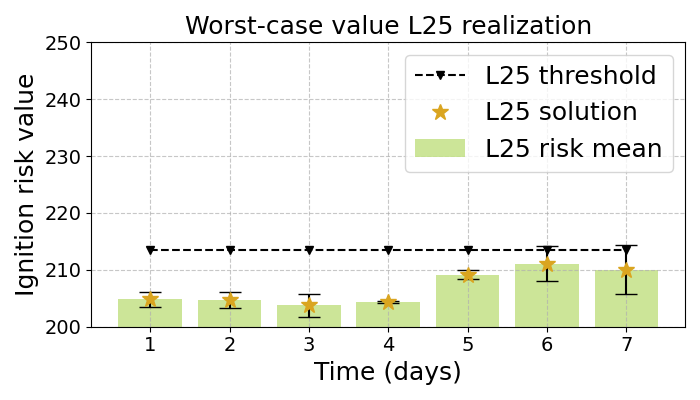}}
    \end{subcaptionbox}
    \caption{Worst-case wildfire ignition risk for two lines under Scheme 1, representative week 2, and two uncertainty budget cases: (Top) Line 5; (Bottom) Line 25.}
    \label{fig:all_cases_scheme_11}
\end{figure}

The wildfire ignition risk of transmission lines that lead to their de-energized influences the investment decisions and changes with the uncertainty budget. In Figure \ref{fig:all_cases_scheme_11}, for Scheme 1 and $\Gamma_\alpha=1$ and $\Gamma_r=1$, we observe that line L5 is de-energized once during the representative week 2, whereas for $\Gamma_r=0.1$ the same line is not affected during the same representative week. This de-energization triggers a new storage investment at node N5 (south-east node in Figure \ref{fig:maps_scheme_2}). Similarly, line L25 is not considered for de-energization for $\Gamma_r = 0.1$, but is de-energized twice during the representative week 2 with $\Gamma_r = 1$.

\onecolumn
\newpage

\twocolumn

\begin{figure}[htbp]
    \centering
    \begin{subcaptionbox}{$\Gamma_{\alpha}=1$, $\Gamma_r=1$.\label{fig:subfig75}}[.48\linewidth]
        {\includegraphics[width=\linewidth]{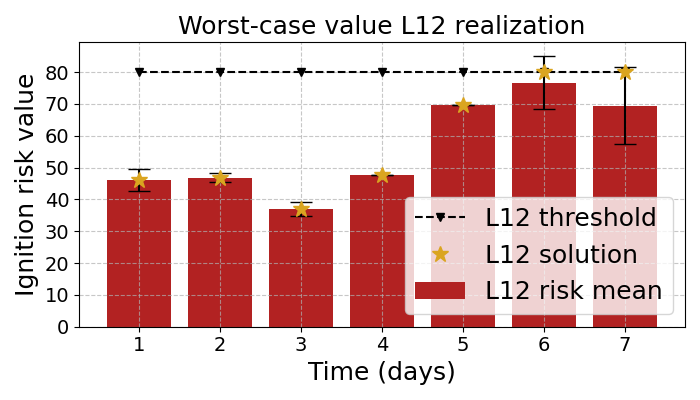}}
    \end{subcaptionbox}
    \hfill
    \begin{subcaptionbox}{$\Gamma_{\alpha}=1$, $\Gamma_r=0.1$.\label{fig:subfig76}}[.48\linewidth]
        {\includegraphics[width=\linewidth]{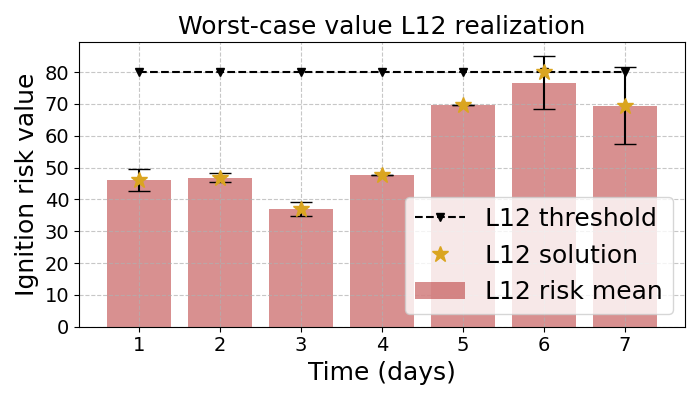}}
    \end{subcaptionbox}
    \vfill
    \begin{subcaptionbox}{$\Gamma_{\alpha}=1$, $\Gamma_r=1$.\label{fig:subfig77}}[.48\linewidth]
        {\includegraphics[width=\linewidth]{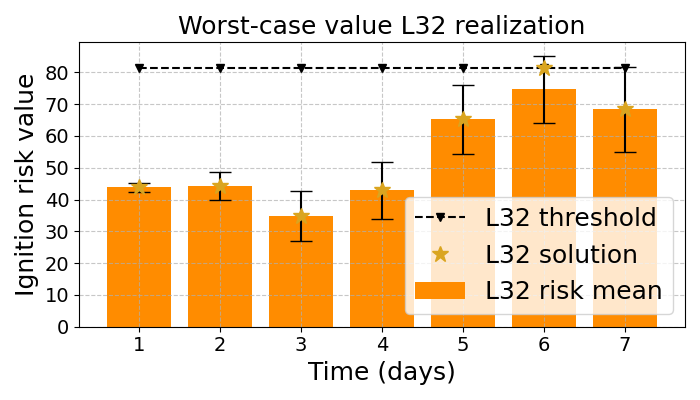}}
    \end{subcaptionbox}
    \hfill
    \begin{subcaptionbox}{$\Gamma_{\alpha}=1$, $\Gamma_r=0.1$. \label{fig:subfig78}}[.48\linewidth]
        {\includegraphics[width=\linewidth]{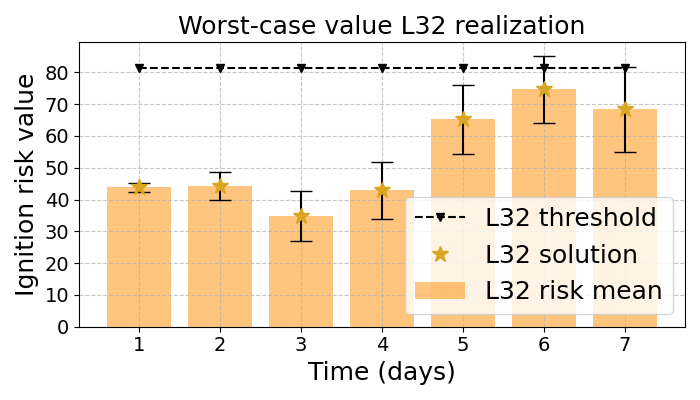}}
    \end{subcaptionbox}
    \caption{Worst-case wildfire ignition risk for two lines under Scheme 1, representative week 2, and two uncertainty budget cases: (Top) Line 12; (Bottom) Line 32.}
    \label{fig:all_cases_scheme_21}
\end{figure}

Figure \ref{fig:all_cases_scheme_21} shows that, under Scheme 2 with $\Gamma_\alpha = 1$ and $\Gamma_r = 1$, line L12 is de-energized once during the representative week 2. However, for $\Gamma_r = 0.1$, the same line is never disconnected. The same behavior is shown for line L32. This consistent de-energization prompts the model to decide to harden this line, aiming to reduce the risk of energy supply interruptions.

\begin{figure}[htbp]
    \centering
    \begin{subcaptionbox}{Frequency of line de-energization with $\Gamma_r=0.1$. \label{fig:subfig92}}[\linewidth]{\includegraphics[width=0.98\linewidth, trim={0, 0, 0, 0}, clip]{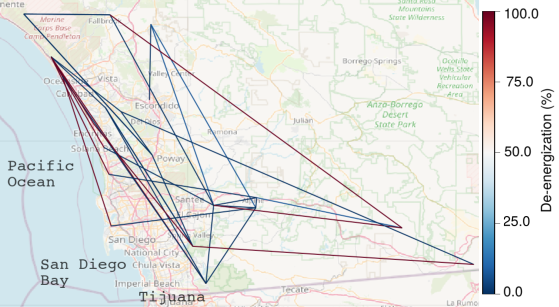}}
    \end{subcaptionbox}
    \begin{subcaptionbox}{Frequency of line de-energization with $\Gamma_r=1.0$.\label{fig:subfig93}}[\linewidth]{\includegraphics[width=0.98\linewidth, trim={0, 0, 0, 0}, clip]{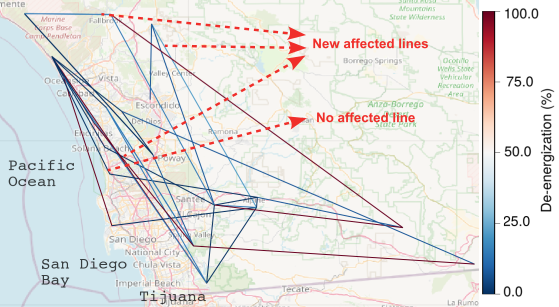}}
    \end{subcaptionbox}
    \caption{Number of times that the line is de-energized during the planning horizon under Scheme 2 and $\Gamma_{\alpha}=1$: (Top) $\Gamma_r=0.1$; (Bottom) $\Gamma_r=1.0$.}
    \label{fig:maps_freq_schemes}
\end{figure}

The summary of de-energization frequency is presented in Figure \ref{fig:maps_freq_schemes}. We observe that the disconnection patterns vary across lines and also change with the selection of uncertainty parameters $\Gamma_{\alpha}$ and $\Gamma_r$, as shown in Figures \ref{fig:all_cases_scheme_11} and \ref{fig:all_cases_scheme_21}. As $\Gamma_r$ increases, de-energization patterns of some lines change. For example, line L23, which is always de-energized when $\Gamma_r=0.1$ is never selected when $\Gamma_r=1$. Instead, lines L5, L11, L15, and L25 (among others) experience one or two de-energization events, all occurring at the same time slots, thereby creating stressful system operating conditions on day 6 of representative week 2.

\begin{figure}[htbp]
    \centering
    \begin{subcaptionbox}{$\Gamma_{\alpha}=0.1$, $\Gamma_r=1$\label{fig:subfig81}}[.48\linewidth]
        {\includegraphics[width=\linewidth, trim={0, 0, 0, 0}, clip]{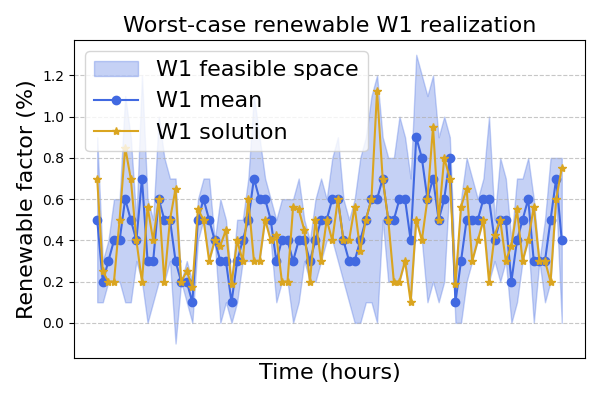}}
    \end{subcaptionbox}
    \hfill
    \begin{subcaptionbox}{$\Gamma_{\alpha}=1$, $\Gamma_r=1$\label{fig:subfig82}}[.48\linewidth]
        {\includegraphics[width=\linewidth, trim={0, 0, 0, 0}, clip]{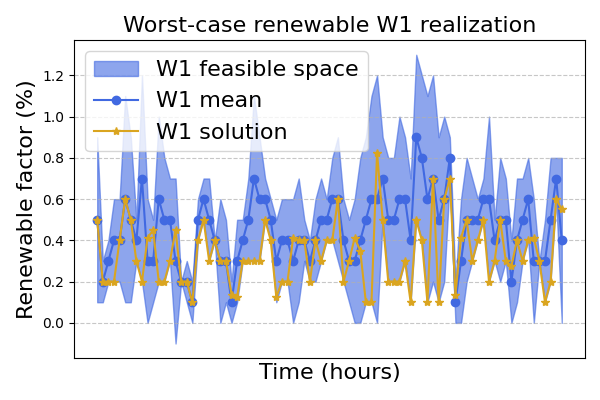}}
    \end{subcaptionbox}
    \vspace{-2.5mm}
    \caption{Worst-case case realization for W1 capacity factor under Scheme 2.}
    \label{fig:representative_21}
\end{figure}

\begin{figure}[htbp]
    \centering
    \begin{subcaptionbox}{$\Gamma_{\alpha}=0.1$, $\Gamma_r=1$\label{fig:subfig83}}[.48\linewidth]
        {\includegraphics[width=\linewidth, trim={0, 0, 0, 0}, clip]{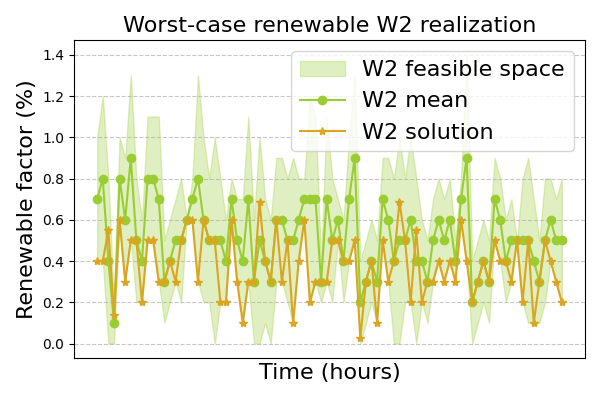}}
    \end{subcaptionbox}
    \hfill
    \begin{subcaptionbox}{$\Gamma_{\alpha}=1$, $\Gamma_r=1$\label{fig:subfig84}}[.48\linewidth]
        {\includegraphics[width=\linewidth, trim={0, 0, 0, 0}, clip]{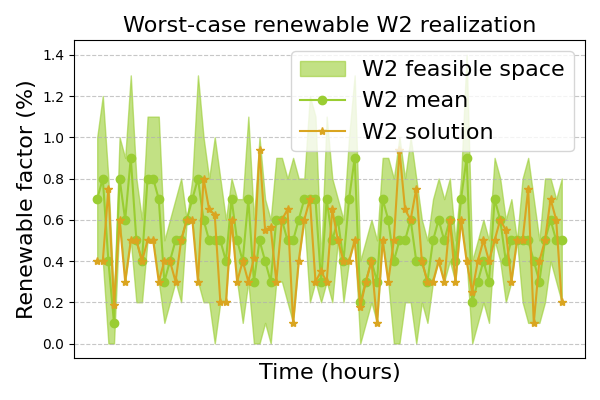}}
    \end{subcaptionbox}
    \vspace{-1.5mm}
    \caption{Worst-case case realization for W2 capacity factor under Scheme 2.}
    \label{fig:representative_22}
\end{figure}

Different values of $\Gamma_{\alpha}$ affect the worst-case realizations of renewable energy generation, as shown in Figures \ref{fig:representative_21} and \ref{fig:representative_22} for renewable generators W1 and W2, respectively. Figure \ref{fig:representative_21} shows that W1 decrease the overall generation, with a greater magnitude during the second representative week. A similar pattern is observed in Figure \ref{fig:representative_22}, where the optimal worst-case capacity-factor of W2 decreases $18\%$ when $\Gamma_{\alpha}=1$.

\section{Conclusion}
This paper develops a practical decision-support model based on optimization for highly renewable power systems. The proposed approach integrates a robust wildfire simulation model that identifies worst-case wildfire risks for line de-energization and the associated worst-case renewable power profiles over a given time horizon.

The model enables optimal investment decisions in undergrounding transmission lines and sizing and locating battery storage, informed by power system operations using real-world wildfire ignition risk and renewable generation data.

By incorporating two distinct uncertainty sets, the model captures the worst-case realizations of wildfire induced line de-energization, variable renewable generation, and their correlations. Furthermore, we propose an evaluation framework through a sensitivity analysis that compares two different investment schemes under varying uncertainty budgets, using a model of the San Diego power system. The results show that bigger uncertainty budgets lead to larger investment needs, requiring a balance between storage and line hardening. When investments are allowed in both undergrounding and storage, the model tends to prioritize undergrounding transmission lines, reflecting the limited investment locations for batteries and the limited number and critical importance of transmission lines for energy supply reliability.

This work opens several directions for future research. One is to scale the model to larger systems, such as the CATS system \cite{taylor2023california}, which represents a 9,000-bus model of the CAISO system. To improve scalability, the proposed algorithm could be enhanced through decomposition and parallel computing, particularly by decomposing the problem by representative week and applying parallel solution procedures such as progressive hedging to coordinate the resulting subproblems while preserving tractability \cite{piansky2024long}. Another direction is to enhance modeling fidelity by incorporating more accurate optimal power flow formulations and additional sources of uncertainty, such as demand variability or other climate-driven natural disasters. Finally, future work could explore integrating a broader range of mitigation actions against wildfire risk and developing computational enhancements for parallelizing the iterative solution and decision process flow.

\begin{acks}
This work was performed under the auspices of the U.S. Department of Energy by Lawrence Livermore National Laboratory under Contract DE-AC52-07NA27344. The work was supported LLNL LDRD Program under Project No. LDRD25-SI-007. The authors gratefully acknowledge Gurobi Optimization for providing a free academic license for their software, which was used in this work.
\end{acks}

\appendix
\section{Appendix: Notation}
\label{sec:nomenclature}

\subsection*{Sets}
\begin{description}
    \item[$G$] Set of thermal generators, indexed by $b$.
    \item[$L$] Set of transmission lines, indexed by $l$.
    \item[$R$] Set of renewable generators, indexed by $r$.
    \item[$S$] Set of battery storage, indexed by $s$.
    \item[$S^C$] Set of candidate battery storage projects, with $S^C\subset S$.
    \item[$B$] Set of buses, indexed by $b$.
    \item[$T$] Set of hours, indexed by $t$.
    \item[$D$] Set of days, indexed by $d$.
    \item[$W$] Set of weeks, indexed by $w$.
\end{description}

\subsection*{Parameters}
\begin{description}
    \item[$p^{max}_g$] Maximum power output (MW).
    \item[$p^{min}_g$] Minimum power output (MW).
    \item[$p^{max}_r$] Maximum renewable power output (MW).
    \item[$P_s$] Maximum charging and discharging rate (MW).
    \item[$\eta$] Battery efficiency (\%).
    \item[$B_l$] Reactance of transmission line $l$ [$\Omega$].
    \item[$C_g$] Linear cost of thermal generator $g$ (\$)  .
    \item[$C_s^{\rm dis}$] Discharge cost of energy storage system $s$ (\$).
    \item[$C^{\rm LS}$] Load shedding cost (\$).
    \item[$\overline{f}_l$] Maximum power flow on line $l$ (MW).
    \item[$X_s^{\rm sup}$] Maximum capacity storage investment (MW).
    \item[$\alpha_r^{t,d,w}$] Availability factor for the renewable generator $r$ (\%).
    \item[$\alpha_{r}^{\text{avg},t,d,w}$] Availability factor mean for the generator $r$ (\%).
    \item[$\alpha_{r}^{\text{dev},t,d,w}$] Availability factor deviation for the generator $r$ (\%).
    \item[$M$] Big-M value.
    \item[$\Gamma_l$] Ignition risk uncertainty set budget.
    \item[$r_l^{\text{avg},d,w}$] Ignition risk mean for the line $l$ (\%).
    \item[$r_l^{\text{dev},d,w}$] Ignition risk deviation for the line $l$ (\%).
    \item[$r_l^{\rm th}$] Ignition risk threshold for line $l$ de-energization (\%).
    \item[$b^{\rm in}(l)$] Bus defined as the origin of transmission line $l$.
    \item[$b^{\rm out}(l)$] Bus defined as the end of transmission line $l$.
\end{description}

\subsection*{Variables}
\begin{description}
    \item[$p_g^{t,d,w}$] Power generation output of generator $g$ (MW).
    \item[$p_b^{t,d,w}$] Load shedding at bus $b$ (MW).
    \item[$c_s^{t,d,w}$] Charge energy of storage $s$ (MW).
    \item[$d_s^{t,d,w}$] Discharge energy of storage $s$ (MW).
    \item[$e_s^{t,d,w}$] State of charge of storage $s$ (MW).
    \item[$f_l^{t,d,w}$] Power flow through line $l$ (MW).
    \item[$\theta_i^{t,d,w}$] Phase angle at bus $b$ (rad).
    \item[$z_l^{\rm UG}$] Binary variable for underground investment on line $l$.
    \item[$X_s^{\rm max}$] Capacity investment for storage $s$ (MW).
    \item[$z_l^{d,w}$] Binary variable for de-energization of line $l$.
    \item[$\gamma_l^{d,w}$] Deviation from ignition risk mean for line $l$ (\%).
    \item[$r_l^{d,w}$] Apparent ignition risk at line $l$ (\%).
    \item[$r_l^{\star d,w}$] Real ignition risk at line $l$ (\%).
\end{description}

\newpage
\bibliographystyle{ACM-Reference-Format}
\bibliography{bib}

\end{document}